\documentclass{style/nseJournal}

\usepackage{amsmath}
\usepackage{amssymb}
\usepackage{bm}
\usepackage{graphicx}
\usepackage[justification=centering]{caption}
\usepackage{subcaption}
\usepackage{multirow}

\newcommand{\ddt}[1]{\dfrac{\partial #1}{\partial t}}
\DeclareMathOperator*{\argmin}{argmin}

\begin{document}

\title{Proper Orthogonal Decomposition Mode Coefficient Interpolation:
	A Non-Intrusive Reduced-Order Model for Parametric Reactor Kinetics} 

\addAuthor{\correspondingAuthor{Zachary K. Hardy}}{a}
\correspondingEmail{zhardy@lanl.gov}
\addAuthor{Jim E. Morel}{b}

\addAffiliation{a}{Los Alamos National Laboratories, X-Theoretical Design Division\\ P.O. Box 1663, Los Alamos, New Mexico 87545}
\addAffiliation{b}{Texas A\&M University, Department of Nuclear Engineering\\ AI Engineering Building\\ 423 Spence Street, MS 3133, College Station, Texas 77843-3133}

\addKeyword{reduced-order modeling}
\addKeyword{proper orthogonal decomposition}
\addKeyword{non-intrusive}
\addKeyword{radiation transport}
\addKeyword{reactor kinetics}
\addKeyword{parametric}
\addKeyword{interpolation}

\titlePage

\begin{abstract}
In this paper, a non-intrusive reduced-order model (ROM) for parametric reactor kinetics simulations is presented.
Time-dependent ROMs are notoriously data intensive and difficult to implement when nonlinear multiphysics phenomena are considered.
These challenges are exacerbated when parametric dependencies are included.
The proper orthogonal decomposition mode coefficient interpolation (POD-MCI) ROM presented in this work can be constructed directly from lower-dimensional quantities of interest (QoIs) and is independent of the underlying model.
This greatly alleviates the data requirement of many existing ROMs and can be used without modification on arbitrarily complex models or experimental data.
The POD-MCI ROM is demonstrated on a number of examples and yields accurate characterizations of the QoIs within the selected parameter spaces at extremely attractive computational speed-up factors relative to the full-order models (FOMs).
\end{abstract}

\section{Introduction}
\label{sec: Introduction} 

Many modern reduced-order models (ROMs) are intrusive, requiring access to, and modification of, the high-dimensional operators of the full-order model (FOM) they seek to emulate.
These techniques often use projection methods to reduce the dimension of the FOM from $N$ to $r \ll N$ with a global reduced-basis obtained from a limited number of FOM executions over a predefined temporal or parametric space.
Such model order reduction (MOR) techniques have been studied in computational fluid dynamics since the late 1960s \cite{lumley1967structure}.
With the introduction of Galerkin proper orthogonal decomposition (GPOD) in the late 1980s \cite{sirovich1987turbulenceI, sirovich1987turbulenceII, sirovich1987turbulenceIII}, these methods became popular across many different fields of computational physics.

In reactor physics, MOR techniques have become widely adopted in the last decade.
Reduced-basis methods been studied for parametric $k$-eigenvalue problems \cite{buchan2013pod,sartori2016multi, german2019reduced,german2019application,behne2022minimally}, reactor kinetics \cite{prill2014semi,sartori2014comparison}, optimizing angular discretizations \cite{soucasse2019angular,sun2020pod, hughes2020discontinuous}, and in parametric fixed source radiation transport applications \cite{behne2022minimally,tano2021affine}.
With advances in multiphysics platforms such as the Multiphysics Object Oriented Simulation Environment (MOOSE) from Idaho National Laboratory (INL) \cite{permann2020moose} and GeN-Foam$^\text{\textregistered}$ \cite{fiorina2015gen}, complex multiphysics reactor models are becoming much more commonplace.
Intrusive ROMs for such models present a number of distinct challenges.

One such challenge is developing a scheme to project physics couplings and nonlinearities to the reduced-order system.
When applying projection methods to nonlinear terms, their evaluation still requires operations that are dependent upon the dimension of the full-order problem.
This can severely diminish any computational savings that would otherwise be obtained by the ROM.
To evaluate these terms on the reduced-order level, techniques such as the discrete empirical interpolation method (DEIM) are often utilized \cite{chaturantabut2010nonlinear}.
This technique projects the nonlinear function onto an low-rank subspace that approximates the space the nonlinearity maps to, then defines an interpolation scheme in the low-dimensional space.
There are several examples of successful multiphysics ROM implementations \cite{sartori2016multi, german2019application, prill2014semi}, however, these implementations are model-dependent and require significant effort to implement.

Perhaps the greatest challenge of intrusive ROMs is their data requirement as problems grow large.
Projection techniques require the reduced-basis elements to be of the same dimension of the FOM.
The reduced-basis is generally computed via a singular value decomposition (SVD) of snapshots of the FOM in parametric space, instances of time, or both.
When nonlinearities are present, snapshots of the full-order nonlinearities must also be stored and their SVD taken.
When considering a realistic parametric multiphysics reactor transient model, it is easy to see how the data requirements to construct the ROM can grow unmanageable.

In most cases, ROMs are used to determine sensitivities of particular quantities of interest (QoIs) or to find parametric configurations that maximize or minimize particular QoIs.
For intrusive methods, the ROM output is generally a reduced-order representation of the full-order solution.
This then requires the QoI to be computed from the ROM-predicted full-order solution.
Non-intrusive MOR techniques offer an attractive alternative in that they are entirely uncoupled from the FOM.
Using available simulation or experimental data, these techniques simply define functional relationships between FOM inputs and outputs.
This provides the distinct advantage of modeling QoIs directly rather than emulating the FOM.

Many non-intrusive ROMs in reactor physics have focused on lower-dimensional quantities such as does mapping \cite{khuwaileh2020gaussian}, material degradation \cite{baraldi2015prognostics}, and inverse uncertainty quantification (UQ) of system parameters with low-dimensional models \cite{wu2016inverse}, and fuel cycle sensitivity analysis.
Methods such as dynamic mode decomposition (DMD) have been applied to high-fidelity reactor simulations \cite{hardy2019dynamic, mcclarren2019calculating, di2020dynamic, abdo2019modeling, abdo2018analysis}, however, the formulation of this method makes its application to parametric problems a challenge.
To the author's knowledge, non-intrusive MOR techniques have been largely unexplored in high-fidelity parametric reactor kinetics simulations.

This paper presents the proper orthogonal decomposition mode coefficient interpolation (POD-MCI) method, which provides a portable, less data intensive alternative to modern intrusive methods.
Whereas intrusive techniques project POD modes onto the FOM to form a system of $r$ equations for the POD mode coefficients, the POD-MCI method formulates a system of $r$ uncoupled interpolation problems.
The method yields accurate characterizations of QoIs across a number of examples at extremely attractive computation speed-up factors.
The following section describes the POD-MCI method in detail.
Following this, the details of the FOM are presented.
Numerical results for two different test cases will then be shown.
Lastly, the results of the numerical examples will be summarized and some concluding remarks about the method made.

\section{POD Mode Coefficient Interpolation Method}
\label{sec: POD-MCI}

Assume a FOM which models the arbitrary scalar or vector function $f = f(\bm{r}, t; \bm{\mu})$ of spatial coordinate $\bm{r}$, time $t$, and model parameters $\bm{\mu}$.
$f$ may be as simple as a scalar temperature or as complex as multiphysics vector.
Now, assume a multi-query application where the behaviors of a particular $\text{QoI} = \text{QoI}(f(\bm{r}, t; \bm{\mu}))$ within some $d$-dimensional parameter space $\bm{\mu} = [\mu_1, \ldots, \mu_d] \in \mathbb{R}^d$ are desired.
The POD-MCI method seeks to emulate the behaviors of $\text{QoI}(f)$ within the parameter space independent of the FOM.

The POD-MCI method is broken into an offline and online stage.
In the offline stage, QoI data from the FOM across the parameter space of interest is obtained and the ROM is constructed.
This is computationally expensive procedure, but is only performed once.
The online stage is computationally cheap and simply involves the querying of the ROM to characterize the parametric behaviors of the QoI.

Let the $f$ with $N_c$ solution components be discretized with $N_r$ spatial nodes $\vec{\bm{r}} = [\bm{r}_1, \ldots, \bm{r}_{N_r}]$ and $N_t + 1$ temporal snapshots at times $t_n = t_0 + n \Delta t$ for $n = 0, \ldots, N_t$, with fixed $\Delta t$.
The discrete solution at time $t_n$ is then given by $\vec{u}(t_n) = \vec{u}_n \in \mathbb{R}^{N_c N_r}$ and the time-dependent solution by $\bm{X} = [\vec{u}_0, \ldots, \vec{u}_{N_t}] \in \mathbb{R}^{N_c N_r \times N_t + 1}$.
For the remainder of this section, the POD-MCI method is described in terms of the full-order solutions $\bm{X}(\bm{\mu})$.
It should be noted that the ROM is formed in an equivalent manner for any QoI.

To obtain the reduced-basis used in the ROM, the FOM must be sampled over the space of interest $\bm{\mu}$.
The number of samples required is largely dependent upon the behaviors of the QoI within the space and the desired accuracy of the global reduced-basis.
This is discussed in more detail later.
In the simplest case, a tensor-product sampling may be used.
For higher-dimensional parameter spaces, techniques such as Latin hypercube sampling \cite{mckay2000comparison} or Smolyak sparse grids \cite{smolyak1963quadrature} may be more appropriate.
Let there be $M$ parametric snapshots for model configurations $\vec{\bm{\mu}} = [\bm{\mu}_1, \ldots, \bm{\mu}_M]$ with simulation results $\bm{Y} = [\bm{X}(\bm{\mu}_1), \ldots, \bm{X}(\bm{\mu}_M)] \in \mathbb{R}^{N_c N_r \times N_t + 1 \times M}$.
It is important to emphasize that the outputs of each simulation must have the same spatial grid and temporal snapshots.
To extract the global reduced basis, each simulation result must be flattened such that $\hat{\bm{X}} = [\vec{u}_0^T, \ldots, \vec{u}_{N_t}^T]^T \in \mathbb{R}^{N_c N_r (N_t + 1)}$.
This results in the matrix
\begin{equation}
	\hat{\bm{Y}} = [\hat{\bm{X}}(\bm{\mu}_1), \ldots, \hat{\bm{X}}(\bm{\mu}_M)] = 
	\begin{bmatrix}
		\vec{u}_0(\bm{\mu}_1) &  \dots & \vec{u}_0(\bm{\mu}_M) \\
		\vdots & \ddots & \vdots \\
		\vec{u}_{N_t}(\bm{\mu}_1) & \dots & \vec{u}_{N_t}(\bm{\mu}_M)
	\end{bmatrix},
	\label{eq: Snapshot Matrix}
\end{equation}
where $\hat{\bm{Y}} \in \mathbb{R}^{N_c N_r (N_t + 1) \times M}$.
Each column of the matrix $\hat{\bm{Y}}$ is referred to as a training snapshot and represents a fixed parametric configuration $\bm{\mu}_i$.
Each row of the matrix, therefore, represents the parametric behaviors of a particular discrete simulation output.
With $N_c N_r (N_t + 1) \gg M$, the matrix $\hat{\bm{Y}}$ is long and skinny.
For realistic simulations, it is easy to see that the dimensionality of a full-order solution can easily become infeasible.
Whereas intrusive methods require data of this dimensionality to be stored, the dimensionality of the matrix $\hat{\bm{Y}}$ can be dramatically reduced with the POD-MCI method by using a lower-dimensional QoI.

To extract the global reduced-basis, the standard POD procedure is employed.
This procedure is essentially a coordinate transformation to a lower-dimensional space that optimally describes the parametric behaviors of the QoI.
The POD basis $\bm{\Phi} = \{ \varphi \}_{i=1}^r$ satisfies the minimization problem
\begin{equation}
	\min_{\bm{\Phi}} \lVert \hat{\bm{Y}} - \mathcal{P}_{\bm{\Phi}} \hat{\bm{Y}} \rVert_2^2 = \min_{\bm{\Phi}} \sum_{i=1}^{M} \int \int \lVert \hat{\bm{X}}(\bm{r}, t; \bm{\mu}_i) - \mathcal{P}_{\bm{\Phi}} \hat{\bm{X}}(\bm{r}, t; \bm{\mu}_i) \rVert_2^2 \, d\bm{r} dt,
\end{equation}
where $\mathcal{P}_\Phi$ is an orthogonal projection onto the subspace spanned by $\bm{\Phi} = \{ \varphi_i \}_{i=1}^r$ and $r \le M$ is the rank of the subspace \cite{rathinam2003new}.
The simulation data $\hat{\bm{Y}}$ is then expressed as $\hat{\bm{Y}}(\bm{r}, t; \bm{\mu}) = \bm{\Phi}(\bm{r}, t) \vec{\bm{a}}(\bm{\mu})$, or
\begin{equation}
	\hat{\bm{X}}(\bm{r}, t; \bm{\mu}_k) = \sum_{i=1}^{r} a_i(\bm{\mu}_k) \varphi_i(\bm{r}, t) = \bm{\Phi}(\bm{r}, t) \bm{a}(\bm{\mu}_k), \hspace{0.25cm} k = 1, \ldots, M,
	\label{eq: POD Expansion}
\end{equation}
where $\bm{\Phi} \in \mathbb{R}^{N_c N_r (N_t + 1) \times r}$ are the POD modes and $\vec{\bm{a}} \in \mathbb{R}^{r \times M}$ are the POD mode coefficients, or coordinates, of the $M$ parametric snapshots.
The basis $\bm{\Phi}$ is obtained via the SVD such that
\begin{equation}
	\hat{\bm{Y}} = \bm{U} \bm{\Sigma} \bm{V}^*
	\label{eq: SVD}
\end{equation}
where $\bm{U} \in \mathbb{R}^{N_c N_r (N_t + 1) \times r}$ are the orthonormal left singular vectors, $\bm{\Sigma} = \text{diag}([\sigma_1, \ldots, \sigma_r]) \in \mathbb{R}^{r \times r}$ are the non-negative singular values ordered greatest to least such that $\sigma_1 > \sigma_2 > \ldots > \sigma_r$,  and $\bm{V} \in \mathbb{R}^{M \times r}$ are the orthonormal right singular vectors.
The operator $*$ denotes the complex conjugate transpose.
The POD modes are the left singular vectors and the singular values a measure of how much ``information'' within $\hat{\bm{Y}}$ each POD mode contains.

Yet to be discussed is the selection of the rank $r$.
There are several common techniques that can be used to select this.
First, a rank $r \in [1, M]$ can simply be selected \emph{a priori}.
Second, a singular value cutoff can be used such that
\begin{equation}
	r = \argmin_{r \in [1, M]} \left( \frac{\sigma_r}{\max(\sigma)} < \tau_\text{cut} \right),
\end{equation}
where $\tau_\text{cut} \in (0, 1)$ is the relative singular value cutoff threshold such that all less than $\tau_\text{cut}$ are truncated.
Lastly, an energy-based truncation is used such that
\begin{equation}
	r = \argmin_{r \in [1, M]} \left(\dfrac{\sum_{i=1}^{r} \sigma_i^2}{\sum_{i=1}^{M} \sigma_i^2} < 1 - \tau \right),
\end{equation}
where $\tau \in (0, 1)$ is the allowable information, or ``energy,'' truncation.
This technique keeps the fewest number of modes which contain the specified fraction of the total energy.

To complete the offline phase, the data within $\hat{\bm{Y}}$ is transformed to POD coordinates.
This is done using Eq. \ref{eq: SVD} such that
\begin{subequations}
	\begin{gather}
		\hat{\bm{Y}} = \bm{U} \bm{\Sigma} \bm{V}^* = \bm{\Phi} \vec{\bm{a}} 
		\label{eq: SVD to POD} \\
		\vec{\bm{a}} = \bm{\Phi}^\dagger \bm{Y} = \bm{\Sigma} \bm{V}^*,
		\label{eq: Dimensionality Reduction}
	\end{gather}
\end{subequations}
where $\dagger$ is the Moore-Penrose pseudoinverse.
To define the interpolation problems, the POD modes, the POD coordinates of each training snapshot, and the training parameters must be stored.
The compression fraction of data achieved by the POD-MCI ROM is then
\begin{equation}
	c = \dfrac{ r \left( N_c N_r (N_t + 1) + M \right) + d M}{ M N_c N_r (N_t + 1) } \approx \dfrac{r}{M},
\end{equation}
where $d \ll N_c N_r (N_t + 1)$ and $r \ll N_c N_r (N_t + 1)$.
This factor is the compression achieved with respect to only the dimensionality reduction obtained using POD.
Significantly more compression relative to full-order solution is achieved when using lower-dimensional QoIs.

From Eq. (\ref{eq: POD Expansion}), using the functional form of the POD coordinates $\bm{a}(\bm{\mu})$, given an unseen configuration $\bm{\mu}_p$, the POD coordinates of the resulting QoI can be obtained via $r$ independent interpolations
\begin{equation}
	a_i(\bm{\mu}_p) = \mathcal{I}\left( \bm{\mu}_p | \bm{a}_i, \vec{\bm{\mu}} \right), \hspace{0.25cm} i = 1, \ldots, r,
\end{equation}
where $\bm{a}_i \in \mathbb{R}^{M}$ is a vector containing the $i$'th POD coordinate for each of the training snapshots.
Because the POD modes are orthonormal, each POD coordinate is uncorrelated, therefore, each interpolation is uncorrelated.
The choice of interpolant $\mathcal{I}$ is arbitrary.
Examples of standard interpolants are multivariate linear, cubic spline, nearest neighbor, and radial basis functions (RBFs).
With interpolated POD coordinates, the predicted QoI for model configuration $\bm{\mu}_p$ is given by
\begin{equation}
	\hat{\bm{X}}(\bm{\mu}_p) = \bm{\Phi} \bm{a}(\bm{\mu}_p) = \sum_{i=1}^{r} a_i(\bm{\mu}_p) \varphi_i,
\end{equation}
which can then be unstacked and further analyzed, if desired

\section{Full-Order Model}
\label{sec: FOM}

In this paper, a multi-group time-dependent neutron diffusion model with delayed neutron precursors is considered.
This is given in Eq. (\ref{eq: FOM})
\begin{subequations}
	\begin{gather}
		\begin{aligned}
			\frac{1}{v_g} &\ddt{\phi_g(\bm{r}, t)} - \nabla \cdot D_g(\bm{r}) \nabla \phi_g(\bm{r}, t) + \Sigma_{t,g}(\bm{r}) 	\phi_g(\bm{r}, t) = \sum_{g'=1}^{G} \Sigma_{s, g' \rightarrow g}(\bm{r}) \phi_g(\bm{r}, t) +\\  &\chi_g^p(\bm{r}) \sum_{g'=1}^{G} \nu_{g'}^p(\bm{r}) \Sigma_{f,g'}(\bm{r}) \phi_g(\bm{r}, t) + \sum_{j=1}^{J} \chi_{g,j}^d (\bm{r}) \lambda_j(\bm{r}) C_j(\bm{r}, t) + Q_q(\bm{r}, t ), \hspace{0.25cm} g = 1, \ldots, G
		\end{aligned}
	\\
	\ddt{C_j(\bm{r}, t)} + \lambda_j(\bm{r}) C_j(\bm{r}, t) = \gamma_j(\bm{r}) \sum_{g=1}^{G} \nu_g^d(\bm{r}) \Sigma_{f,g} \phi_g(\bm{r}, t), \hspace{0.25cm} j = 1, \ldots J,
	\end{gather}
	\label{eq: FOM}
\end{subequations}
where it is assumed appropriate boundary and initial data are supplied.
Here, $\phi_g(\bm{r}, t)$ is the group $g$ scalar flux and $C_j(\bm{r}, t)$ is the delayed neutron precursor species $j$ concentration at position $\bm{r} \in \bm{\Omega}$ at time $t \in [0, T]$.
The total cross-section, fission cross-section, diffusion coefficient, neutron speed, and external source for group $g$ are given by $\Sigma_{t,g}(\bm{r})$, $\Sigma_{f,g}(\bm{r})$, $D_g(\bm{r})$, $v_g$, and $Q_g(\bm{r}, t)$ respectively.
The differential scattering cross-section from group $g'$ to group $g$ is given by $\Sigma_{s, g' \rightarrow g}(\bm{r})$.
The prompt fission spectrum is given by $\chi_g^p(\bm{r})$, where $\sum_g \chi_g^p(\bm{r}) = 1$, and the delayed emission spectra by $\chi_{g,j}^d(\bm{r})$, where $\sum_g \chi_{g,j}^d(\bm{r}) = 1$, for $j = 1, \ldots J$.
Further, the prompt and delayed neutrons per fission are given by $\nu_g^p(\bm{r})$ and $\nu_g^d(\bm{r})$, respectively.
Lastly, the delayed neutron precursor decay constants are given by $\lambda_j(\bm{r})$ and fission yields by $\gamma_j(\bm{r})$ where $\sum_j \gamma_j(\bm{r}) = 1$.
It should be emphasized that the cross-sections are not restricted to the functional dependence listed above.
In fact, cross-sections may often be modeled as a function of time to encapsulate the behaviors of control rod withdrawals, a function of temperature to include the effects of Doppler broadening, or both.

In this paper, a finite volume (FV) discretization is used to discretize space.
The spatial domain is first divided into $N$ cells.
The multi-group scalar flux is then approximated with cell-wise averages $\phi_{g,i}(t)$.
Integrating over the volume of each cell then gives
\begin{equation}
	\frac{V_i}{v_g} \ddt{\phi_{g,i}(t)} + \sum_{f \in \partial \Omega_i} A_f \bm{J}_{g,f}(t) \cdot \hat{n}_f + \Sigma_{t,g,i} \phi_{g,i}(t) V_i = S_{g,i} V_i,
\end{equation}
where $V_i$ is the volume of cell $i$, $A_f$ is the area of face $f$ of cell $i$ with outward normal $\hat{n}_f$, $S_{g,i}$ is the scattering, fission, and inhomogeneous source, and
\begin{equation}
	\bm{J}_{g,f}(t) = - D_{g} \nabla \phi_g(t) \Big\vert_f
\end{equation}
is the current across face $f$ of cell $i$.
Particular care must be taken in the evaluation of the current.
On orthogonal grids, the gradient dotted with the outward pointing normal reduces to a single derivative, which greatly simplifies matters.
The derivative is approximated with a central difference with respect to the face such that
\begin{equation}
	\nabla \phi(t) \Big\vert_f = \dfrac{\phi_n - \phi_i}{d_{in}},
\end{equation}
where $n$ denotes the neighbor cell and $d_{in}$ denotes the distance between the centers of cell $i$ and its neighbor $n$ with respect to face $f$.
To account for material interfaces, the diffusion coefficient at a cell face is computed using a weighted harmonic mean such that
\begin{equation}
	D_f = \left( \dfrac{w}{D_i} + \dfrac{1 - w}{D_n}  \right)^{-1},
\end{equation}
where $w = d_{if}/d_{in}$, with $d_{if}$ being the distance from the center of cell $i$ to the center of face $f$.
The final discrete system for interior cells is then given by
\begin{equation}
	\frac{V_i}{v_g} \ddt{\phi_{g,i}(t)} + \sum_{f \in \partial \Omega_i} A_f \frac{D_{g,f}}{d_{in}} \left( \phi_{g,i}(t) - \phi_{g,n}(t) \right) + \Sigma_{t,g,i} \phi_{g,i}(t) = S_{g,i}(t) V_i.
\end{equation}
This method is second order and conservative.

The trapezoidal second-order backwards difference formula (TBDF-2) method is used for the temporal discretization \cite{edwards2011nonlinear}.
This method is second-order accurate with strongly damped oscillatory behaviors.
The method begins by taking a half time-step using the Crank-Nicholson method such that
\begin{equation}
	\frac{4}{\Delta t} \left( f^{n+1/4} - f^n \right) = \bm{A} f^{n+1/4}.
\end{equation}
This is then used to compute $f^{n+1/2} = 2 f^{n+1/4} - f^n$.
With this result, a half-time step of second-order backwards difference is taken to obtain the end-of-time step result via
\begin{equation}
	\frac{3}{\Delta t} \left( f^{n + 1} - f^{n+1/2} \right) - \frac{1}{\Delta t} \left( f^{n+1/2} - f^n \right) = \bm{A} f^{n+1}.
\end{equation}

Because the delayed neutron precursor equations are cell-wise ODEs, there is no compelling reason to increase the dimensionality of the multi-group system and solve for them simultaneously.
Instead, inspection of the time-discretized system,
\begin{equation}
	\left( 1 + \lambda_{j,i} \Delta t_* \right)C_{j,i}^* = C_{j,i}^\text{prev} + \gamma{j,i} \Delta t_* \sum_{g=1}^{G} \nu_{g,i}^d \Sigma_{f,g,i} \phi_{g,i}^*,
\end{equation} 
shows that $C_{j,i}^* = C_{j,i}^*(C_{j,i}^\text{prev}, \phi_{g,i}^*)$.
In the above equation, $*$ denotes the unknown time step index, $\Delta t_*$ the effective time step size for a fully implicit version of the method, and the ``prev'' superscript the term involving all previously computed solutions.
For example, for the first half time step of Crank Nicholson $\Delta t_* = \Delta t / 4$. 
With this, $C_{j,i}^*$ can be eliminated from the multi-group diffusion equations and the delayed neutron precursors only need be updated after the linear system is solved.

\section{Numerical Results}
\label{sec: Results}

This section presents numerical results for two examples.
The first is a relatively simple problem and focuses on demonstrating the POD-MCI ROM in detail.
The second is a more realistic application of a reactor pulse.
This is a two-dimensional problem including transient cross-sections and nonlinear temperature feedback with a coupled adiabatic heat-up model.

\subsection{Homogeneous Three-Group Sphere Problem}

This numerical example is a crude model of a pulsed neutron experiment in which a subcritical sphere of fissile material with radius $r_b = 6$ cm is subjected to a burst of fast neutrons.
Because the time scale of such experiments is generally very short, the effects of delayed neutron precursors are omitted.
For a more detailed description of the relevance and motivation behind this problem the reader is referred to \cite{hardy2019dynamic} where this is discussed in detail.

The three-group nuclear data used for this problem is given in Tables \ref{tab: Nuclear Data} and \ref{tab: Scattering Matrix}.
As the data suggests, all fission events yield neutrons in the fast group and only down-scattering from the fast to the epithermal group is permitted.
\begin{table}[!h]
	\centering
	\caption{The nuclear data for the three-group homogeneous sphere problem.}
	\begin{tabular}{| c | c | c | c | c |}
		\hline
		\textbf{Group} & $\sigma_t$ (b) & $\nu\sigma_f$ (b) & $\chi$ & $v$ (cm/$\mu$s) \\ \hline
		0 & 7.71 & 5.4    & 1 & 2000 \\ \hline
		1 & 50.0 & 60.8 & 0 & 100  \\ \hline
		2 & 25.6 & 28.0 & 0 & 2.2 \\ \hline
	\end{tabular}
	\label{tab: Nuclear Data}
\end{table}
\begin{table}[!h]
	\centering
	\caption{The scattering matrix for the three-group homogeneous sphere problem in barns where each row represents the final group $g$ and each column the initial $g'$.}
	\begin{tabular}{| c || c | c |  c |}
			\hline
			\textbf{Group} & 0 & 1 & 2 \\ \hline\hline
			0 & 4.4 & 0 & 0 \\ \hline
			1 & 1.46 & 13.8 & 0 \\ \hline
			2 & 0 & 0 & 12.0 \\
			\hline
	\end{tabular}
	\label{tab: Scattering Matrix}
\end{table}
This data is formulated such that there are no production terms for the thermal group.
In other words, if no thermal neutrons are present in the initial condition, none will be observed in the problem.
A nominal atom density of $\rho = 0.05$ atoms/b-cm is used to compute the macroscopic cross-sections used in the simulation and the diffusion coefficient is computed via $D_g = (3 \rho \sigma_{t,g})^{-1}$. 

\begin{figure}[!h]
	\centering
	\includegraphics[width=0.7\linewidth]{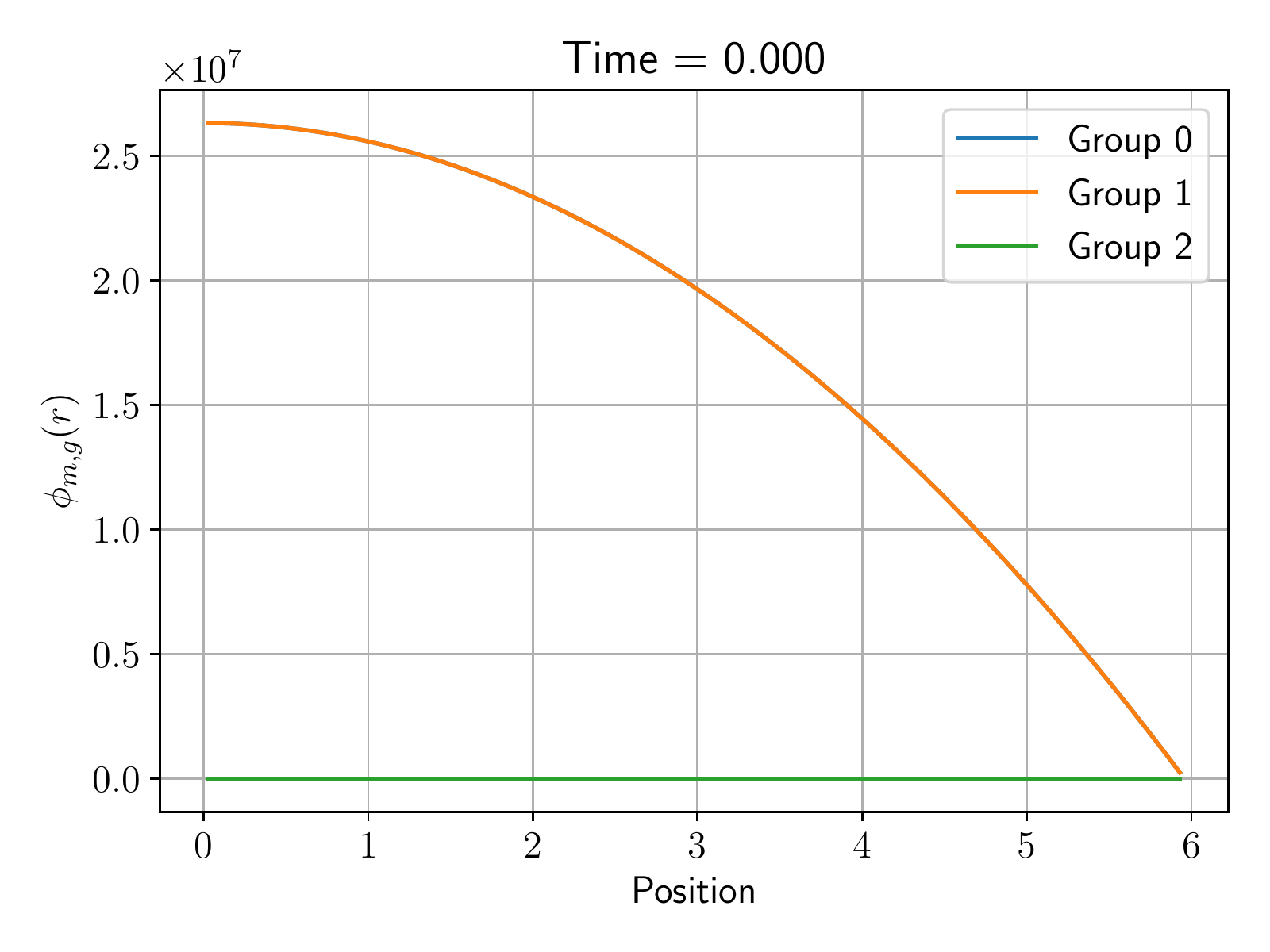}
	\caption{The initial condition for the three-group homogeneous sphere problem.}
	\label{fig: Sphere Initial Condition}
\end{figure}

\begin{figure}
	\centering
	\begin{subfigure}{0.49\linewidth}
		\centering
		\includegraphics[width=\textwidth]{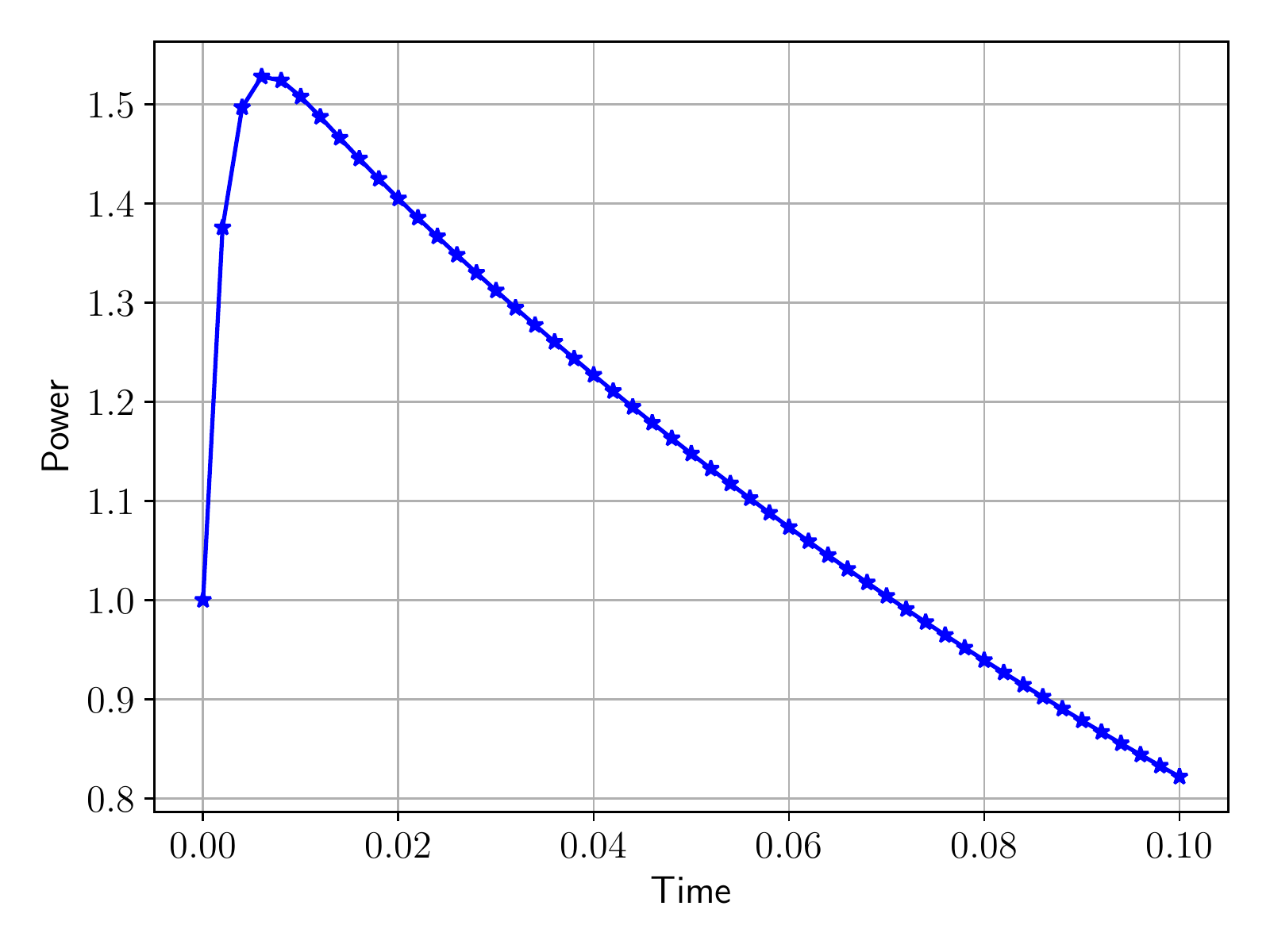}
	\end{subfigure}
	\hfill
	\begin{subfigure}{0.49\linewidth}
		\centering
		\includegraphics[width=\textwidth]{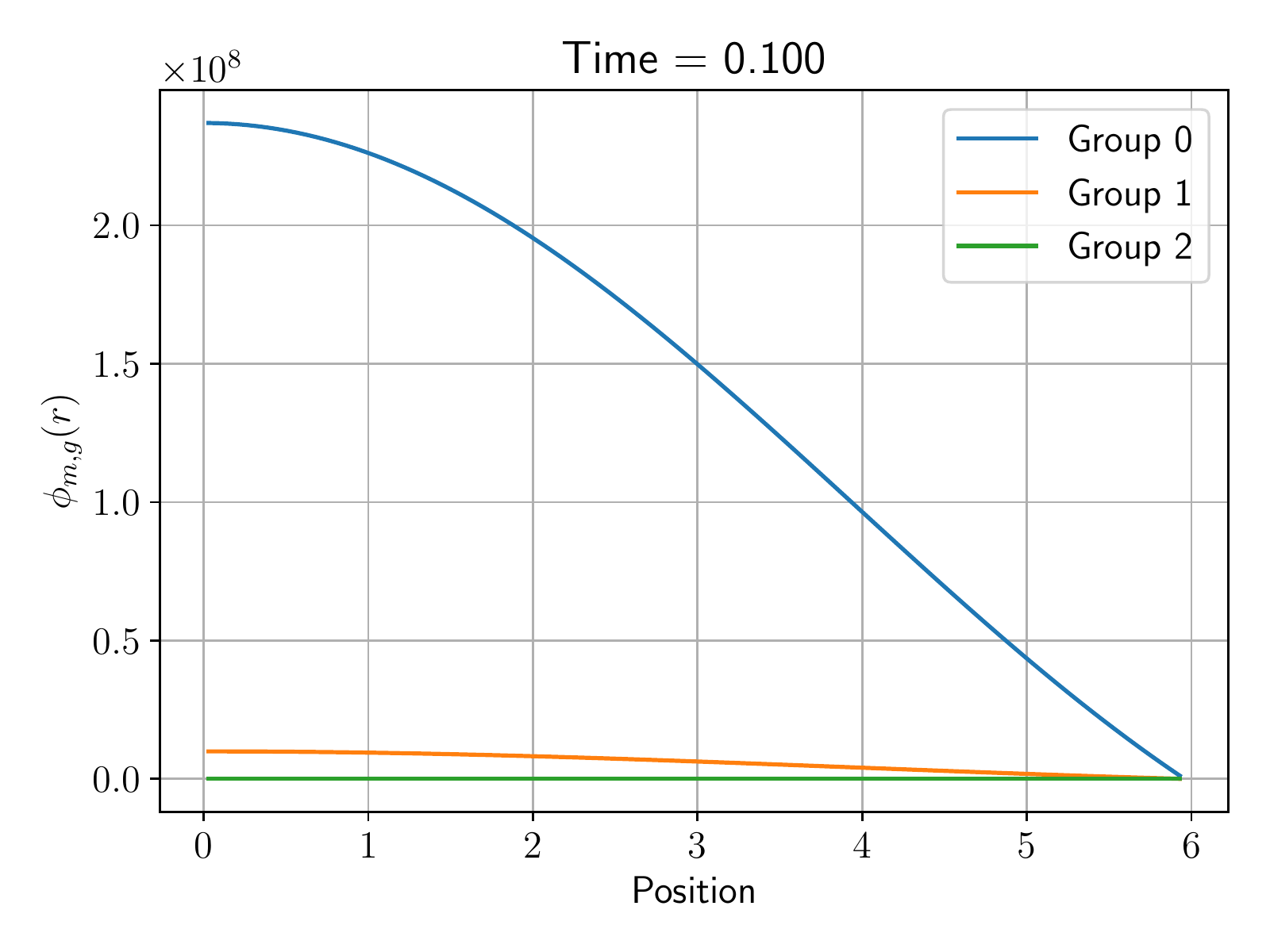}
	\end{subfigure}
	\caption{The power as a function of time (left) and the end of simulation multi-group scalar flux (right).}
	\label{fig: Sphere Power Profile}
\end{figure}

In order to induce more interesting dynamics, a parabolic initial condition in the fast and epithermal groups is used such that
\begin{equation}
	\phi_g(r, t=0) = \begin{cases}
		C \left( 1.0 - \dfrac{r^2}{r^2_b} \right), & \text{if } g = 0, 1 \\
		0.0, & \text{otherwise},
	\end{cases}
\end{equation}
where $C$ is a constant that normalizes the flux to unit power.
This can be seen in Figure \ref{fig: Sphere Initial Condition}.
On the sphere boundary, a zero flux condition is imposed for each group, the simulation is run for 100 ns with time steps of 2 ns, and 100 equally spaces computational cells are used.
With a higher fission cross-section in the epithermal group, a significant number of fission events occur at the start of the problem leading to a spike in the reactor power before settling into its subcritical dynamics.
The reactor power as a function of time and end of simulation scalar flux profile is shown in Figure \ref{fig: Sphere Power Profile}.

The remainder of this section is broken into two subsections.
The first uses a one-dimensional parameterization of the problem to aid in building an intuition about the POD-MCI ROM performance.
In the second section, a more realistic three-dimensional parameterization is considered.

\subsubsection{Parameterization in Sphere Radius}
\label{sec: Sphere 1D}

\begin{figure}[!h]
	\centering
	\includegraphics[width=0.7\linewidth]{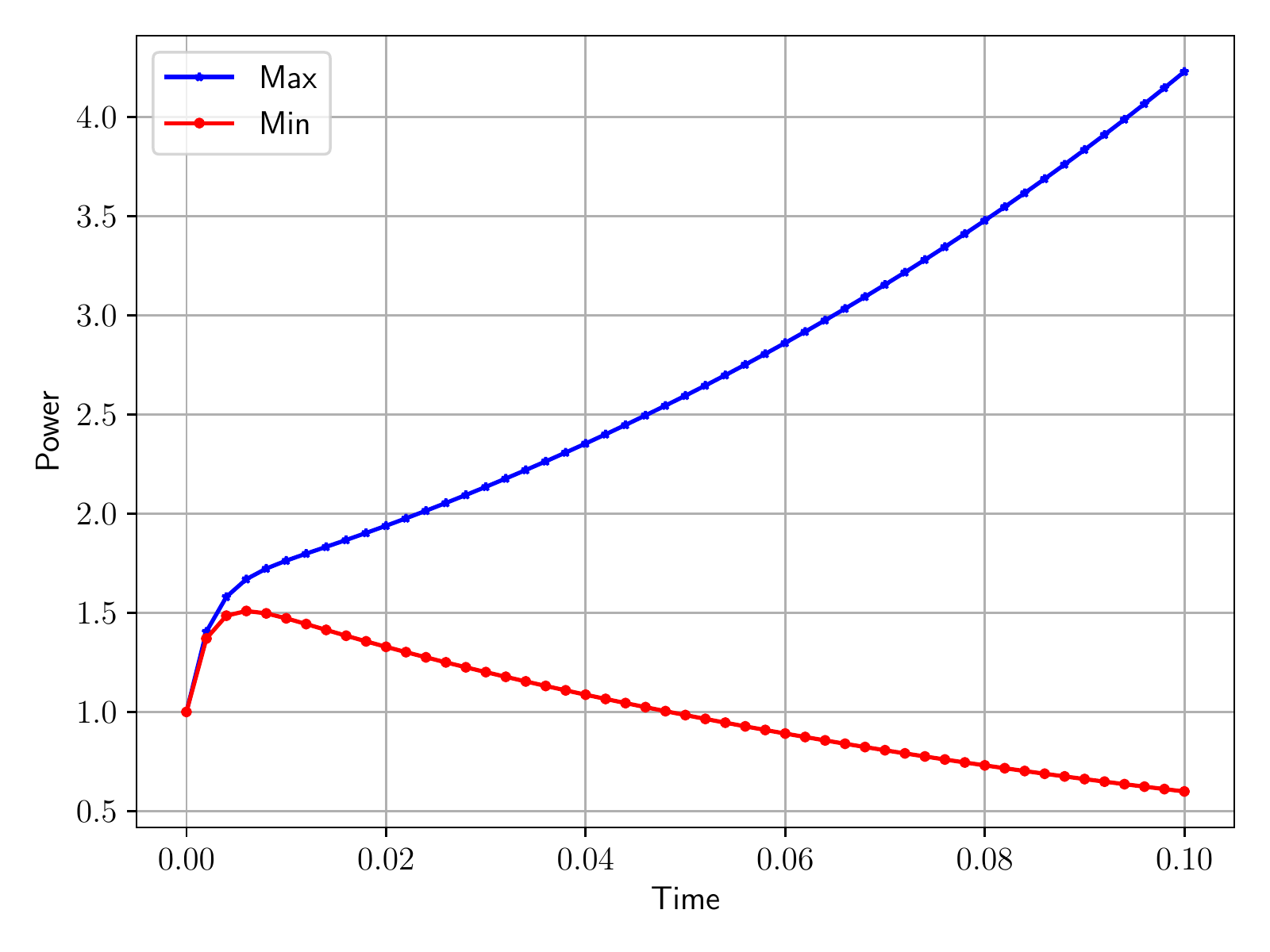}
	\caption{The most subcritical (red) and most supercritical (blue) states in the parameter space.}
	\label{fig: Sphere 1D Power Span}
\end{figure}

In this section the sphere radius is taken to be as the uncertain or design parameter.
In lieu of using the subcritical nominal configuration, the nominal sphere radius is set to the critical radius $r_b = 6.1612$ cm.
It should be highlighted that perturbations about criticality yield exponentially divergent results from the nominal problem.
For this example, 21 simulations are run with radii uniformly spaced between $r_b = 6.1512$ cm $\pm 2.5\%$.
The most subcritical ($r_b - 2.5\%$) and most supercritical ($r_b + 2.5\%$) power profiles are shown in Figure \ref{fig: Sphere 1D Power Span}.
This represents a a 150\% relative difference and a 606\% minimum-to-maximum difference between the final powers.
Further, this spans two vastly different dynamic regimes.

\begin{figure}[!h]
	\centering
	\includegraphics[width=0.7\linewidth]{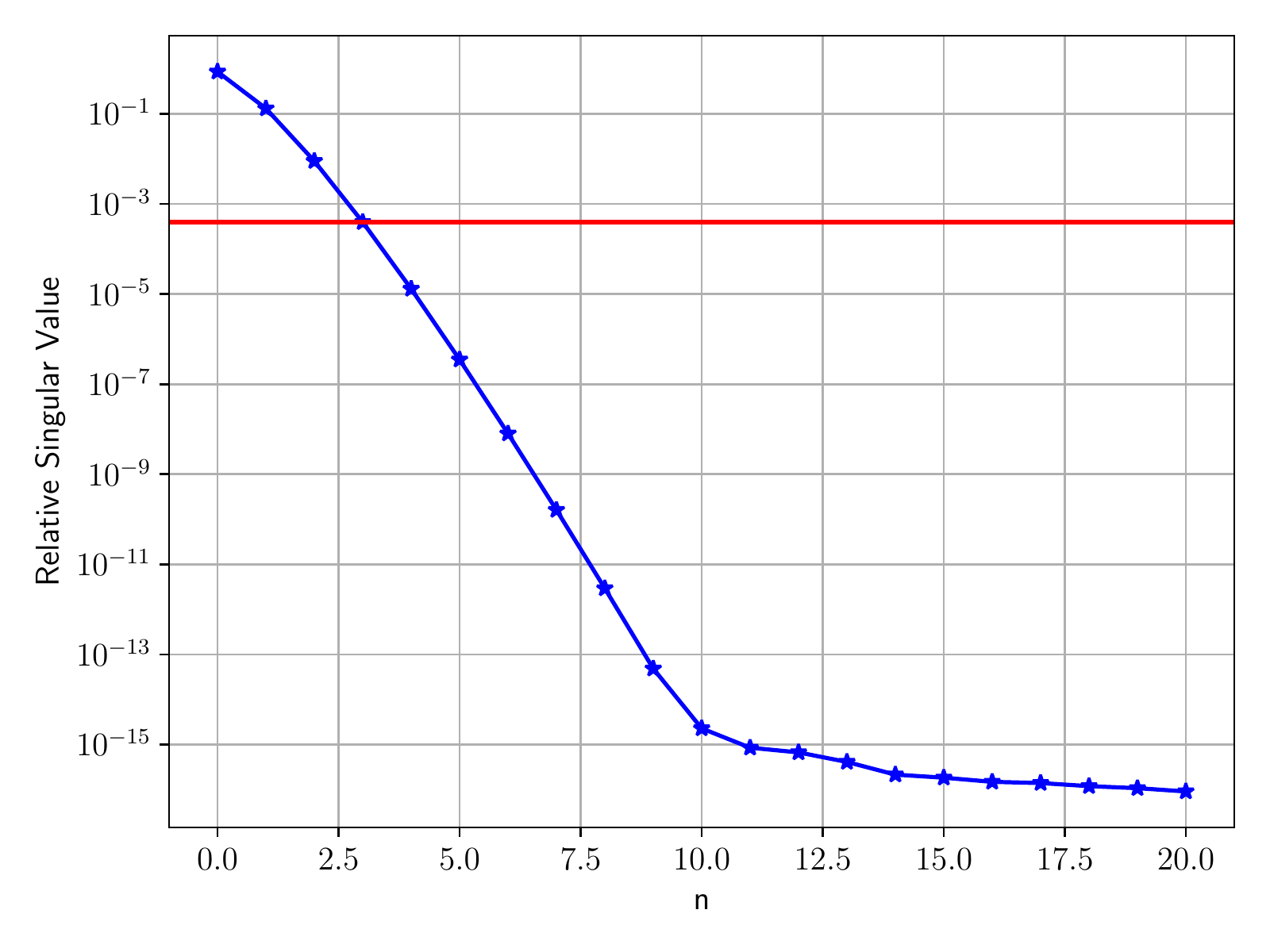}
	\caption{The scree plot for the multi-group scalar flux data.}
	\label{fig: Sphere 1D Scree}
\end{figure}

\begin{figure}[!h]
	\centering
	\includegraphics[width=0.7\linewidth]{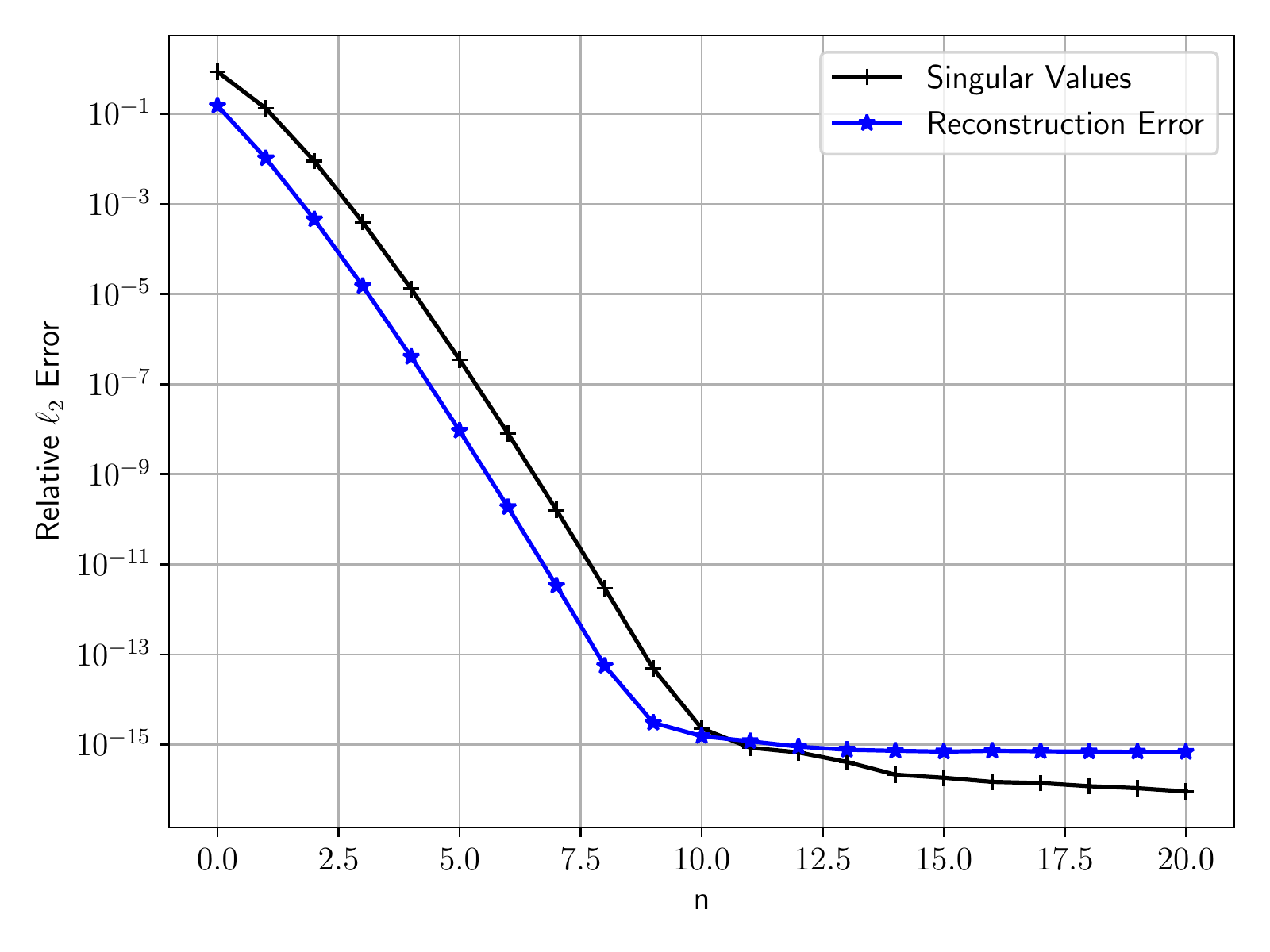}
	\caption{The relative $\ell_2$ POD reconstruction error as a function of mode number plotted against the singular value spectrum.}
	\label{fig: Sphere 1D POD Verification}
\end{figure}

From this data, the POD-MCI ROM is constructed from the multi-group scalar flux data.
The size of each simulation snapshot is the product of the number of time steps, number of grid points, and number of groups, or $51 \times 100 \times 3 = 15,300$.
The energy truncation limit for the ROM is set to $\tau = 10^{-8}$, which leads to a four mode model.
The relative singular value decay, or scree plot, and truncation threshold is shown in Figure \ref{fig: Sphere 1D Scree}.
The relative $\ell_2$ reconstruction error of the training data is $1.5 \times 10^{-5}$, the mean simulation reconstruction error is $1.47 \times 10^{-5}$ and the maximum $4.67 \times 10^{-5}$.
For verification purposes, it is common to investigate the reconstruction error as a function of mode number or energy truncation parameter.
According to theory, the reconstruction error is proportional to the largest truncated singular value such that $\lVert X_r - X \rVert_{\ell_2} \propto \sigma_{r+1}$.
This behavior is observed as shown in Figure \ref{fig: Sphere 1D POD Verification}.
The reconstruction error, of course, is not representative of the POD-MCI error since the mode coefficients are not interpolated, but instead obtained via a least-squares fit to the POD modes.
It does, however, provide an approximate best case scenario.

\begin{figure}[!h]
	\centering
	\includegraphics[width=0.8\linewidth]{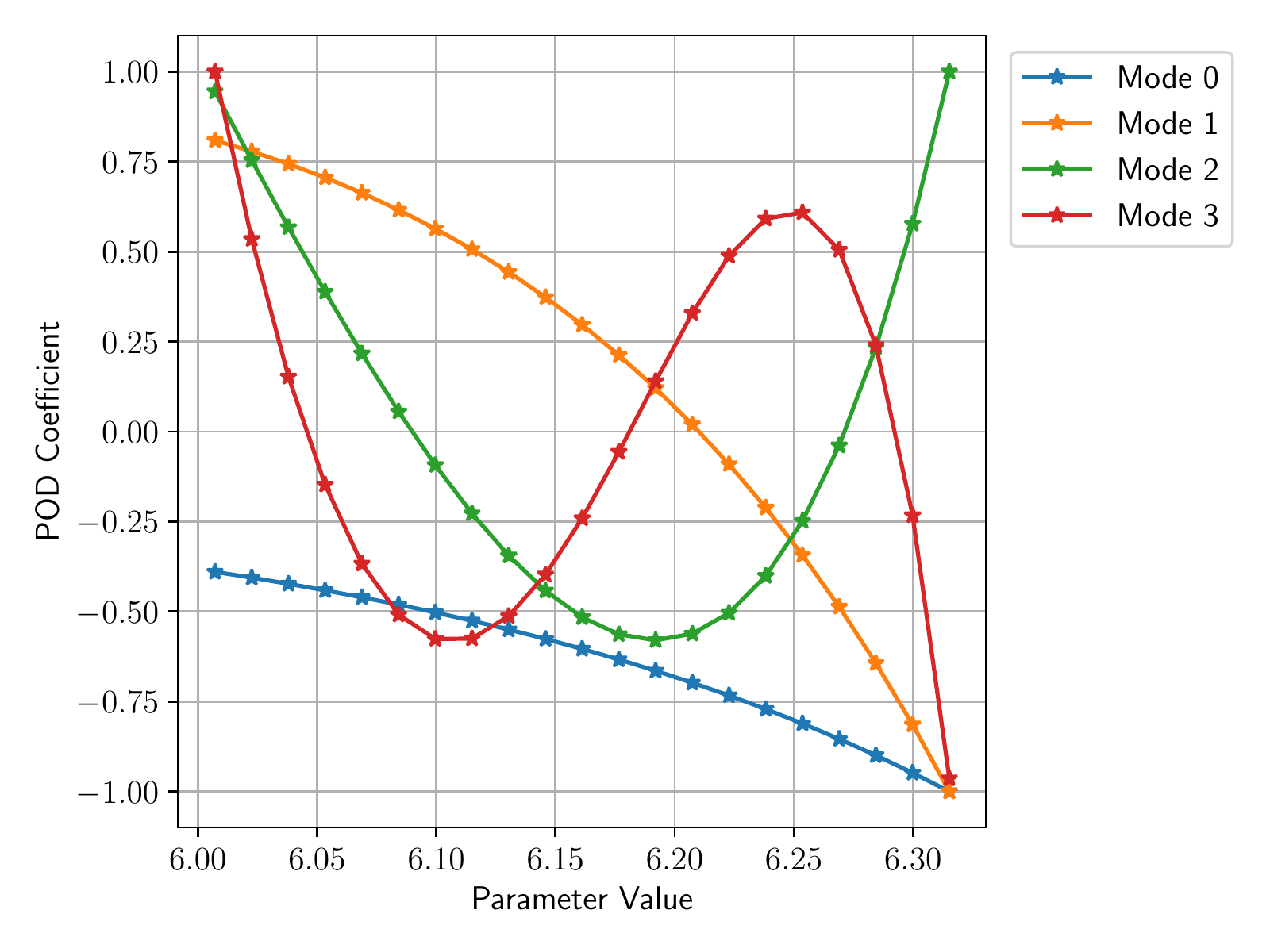}
	\caption{POD mode coefficients as a function of sphere radius.}
	\label{fig: Sphere 1D POD Coeffs}
\end{figure}

It is useful to visualize the POD mode coefficients as a function of parameter value to get an idea of the smoothness of their functional form.
This is shown in Figure  \ref{fig: Sphere 1D POD Coeffs}.
It is clear that the smoothness of the POD modes coefficients degrades with mode number.
While interpolation accuracy will clearly degrade with mode number, due to the properties of POD, this is not entirely troublesome.
Because the POD modes are orthonormal, each successive POD mode adds new, uncorrelated information to the expansion.
Further, because each mode contains successively less information than those prior, the impacts of poor interpolation on higher order modes to the prediction error is significantly damped.

Because the choice of interpolant is arbitrary for the POD-MCI ROM, only radial basis function (RBF) interpolants with a thin-plate spline kernel given by $r^2 \log(r)$ are used in this work.
RBFs are advantageous because no special treatment is required when predicting points outside of the available grid of points.
While extrapolation should always be minimized, in higher dimensional spaces where sampling on regular grids is less feasible, some extrapolation in particular dimensions is unavoidable.
This particular kernel function is used because it does not require the tuning of a shape parameter.
No efforts are made to tune other hyper-parameters.
For optimal results, one should perform a preliminary analysis to determine the optimal interpolant configuration for the particular data set.
In this work, the interpolations are carried out with the \verb*|RBFInterpolator| class within the Python package \verb*|SciPy| \cite{SciPy}.

\begin{figure}[!h]
	\centering
	\includegraphics[width=0.7\linewidth]{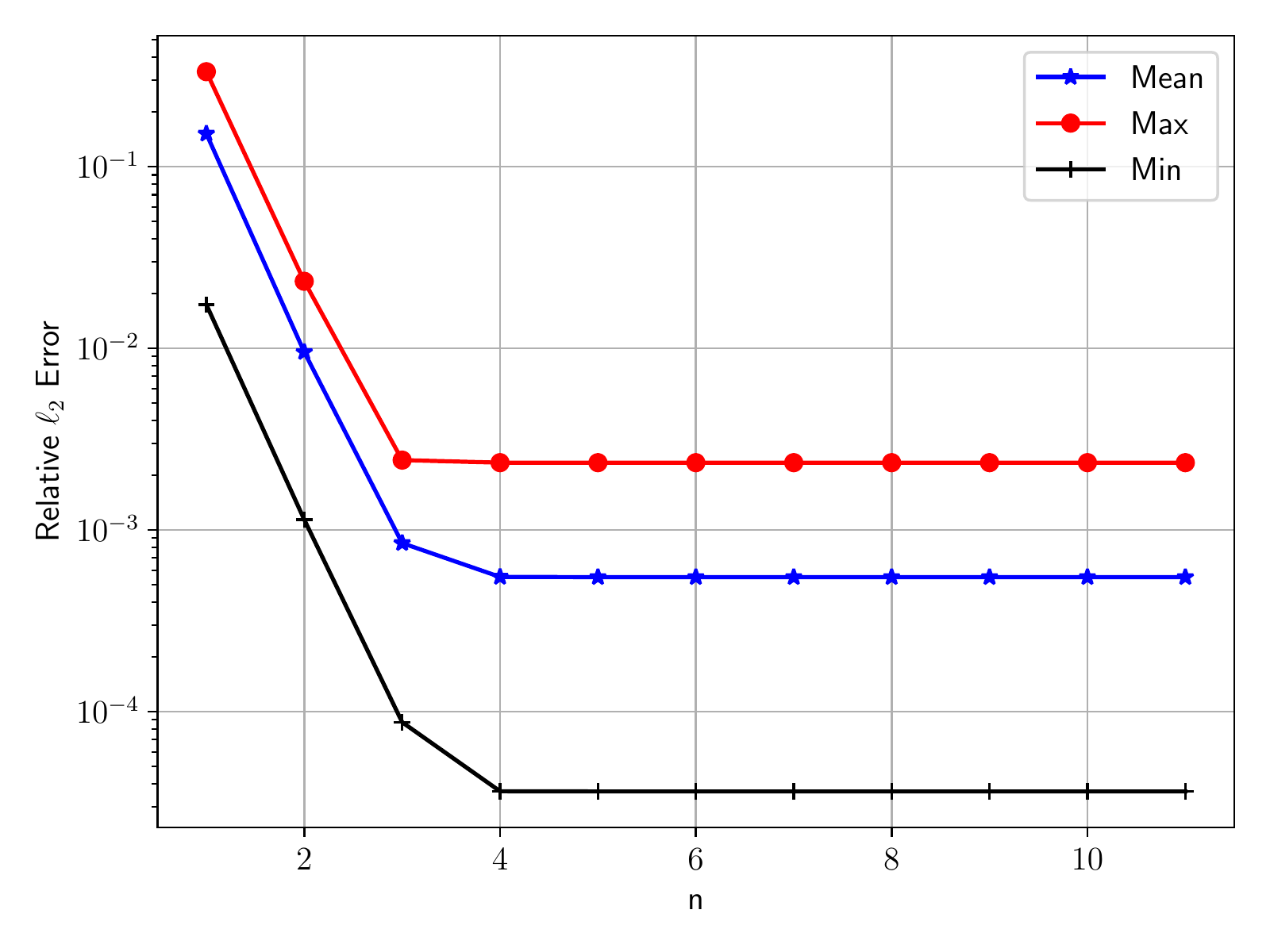}
	\caption{The mean, maximum, and minimum prediction errors as a function of number of POD modes.}
	\label{fig: Sphere 1D POD Trucation Error}
\end{figure}

Taking every other snapshot to beginning with the first to be training data and the remainder to be validation, or test, data, the mean, maximum, and minimum simulation reconstruction error as a function of number of modes is shown in Figure \ref{fig: Sphere 1D POD Trucation Error}.
There are several noteworthy observations from this result.
First, the monotonic decrease in prediction error is demonstrative of the prior statement regarding interpolations on non-smooth higher order POD modes.
Second, for few mode models, the error in the POD representation clearly dominates whereas for several models with more modes, it is clear that the interpolation error dominates.
Further, on a regular parametric grid, using every other snapshot as an interpolation point implies that all interpolations are a distance $\Delta r_b$ away from the nearest training point.
Because interpolation accuracy falls off with distance from training points, this represents a worst case result for grids structures in this way.
It is important to emphasize that these results are obtained with only 11 simulations in the training set.
With more simulation data or an interpolant optimized for this data set, these errors can be improved.

\begin{figure}[!h]
	\centering
	\includegraphics[width=0.8\linewidth]{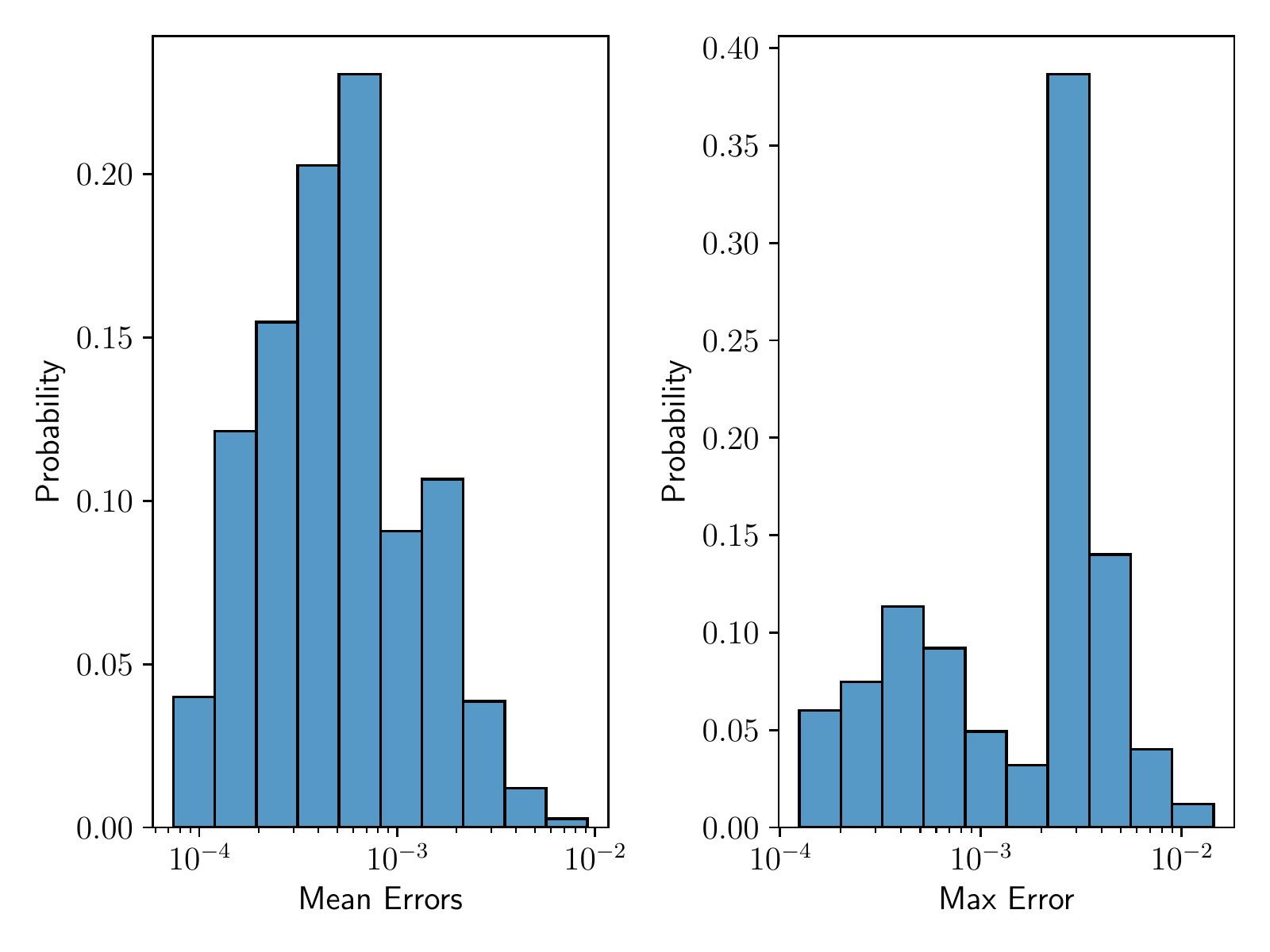}
	\caption{The distribution of mean and maximum errors across 750 cross-validation sets.}
	\label{fig: Sphere 1D Cross Validation}
\end{figure}

To approximate the performance of the POD-MCI ROM on irregular parametric grids, repeated $k$-fold cross-validation is utilized.
This cross-validation technique shuffles the available snapshots and partitions them into $k$ sets, or folds.
Then, $k-1$ of the sets are used to train the ROM while the remaining set is held back for validation.
This procedure is repeated $k$ times where each time a different set is taken as validation.
When only limited data are available, this process can be repeated a number of times, shuffling the snapshots each time to produce different sets.
To carry out this procedure, the tools within the Python package \verb*|Scikit-Learn| are utilized \cite{Scikit-Learn}.
With three sets and 250 repeats for a total of 750 cross-validation sets, the distributions of mean and maximum prediction errors across the ROMs is shown in Figure \ref{fig: Sphere 1D Cross Validation}.
From this, it is clear that regardless of the shuffling, mean ROM errors remain less than 1\% and in a small fraction of cases maximum errors exceed 1\%.
In fact, maximum errors in excess of 1\% lie occurred less than 2.5\% of the time.
When randomly shuffling the training samples, the worst case scenario for maximum error would occur when the training set includes many contiguous samples leaving a significant portion of the parameter space sparsely sampled.
In this case, for a training set with 14 snapshots and a validation set with 7 snapshots, interpolations or extrapolations at a distance of $7 \Delta r_b$ would be theoretically possible.
In practice, with appropriate sampling these scenarios are easily avoided.
For this reason, these distributions are likely conservative estimates of what would be observed in practice.

In real applications, one seeks to use as many training samples as are available to ensure accurate results.
Leave-one-out cross-validation is a useful tool to characterize the performance of the ROM when using all but one set of simulation data.
This methodology is equivalent to $k$-fold cross-validation when $k$ is the number of available snapshots.
While this is limited to only as many cross-validation sets as available snapshots, it gives a good indication of how a ROM will perform given that no large extrapolations from the existing data are explored.
For this data set, the average prediction error was 0.15\% with a maximum of 1.4\% and a minimum of 0.0056\%.

These results demonstrate the ability of the POD-MCI ROM to accurately predict full-order simulation results across vastly different dynamic regimes in a simple one-dimensional parameterization.
Acceptable errors were achieved and found to be generally insensitive to the sampling of the space.
For this problem, each full-order simulation took on average 1.1 s using Python.
For the leave-one-out cross-validation study, the average POD-MCI ROM construction time was 8.12 ms and the average query time 0.135 ms.
This results in a speed-up of approximately 8,000 times which is extraordinarily attractive.

\subsubsection{Three-Dimensional Parameterization}

\begin{figure}[!h]
	\centering
	\includegraphics[width=0.7\linewidth]{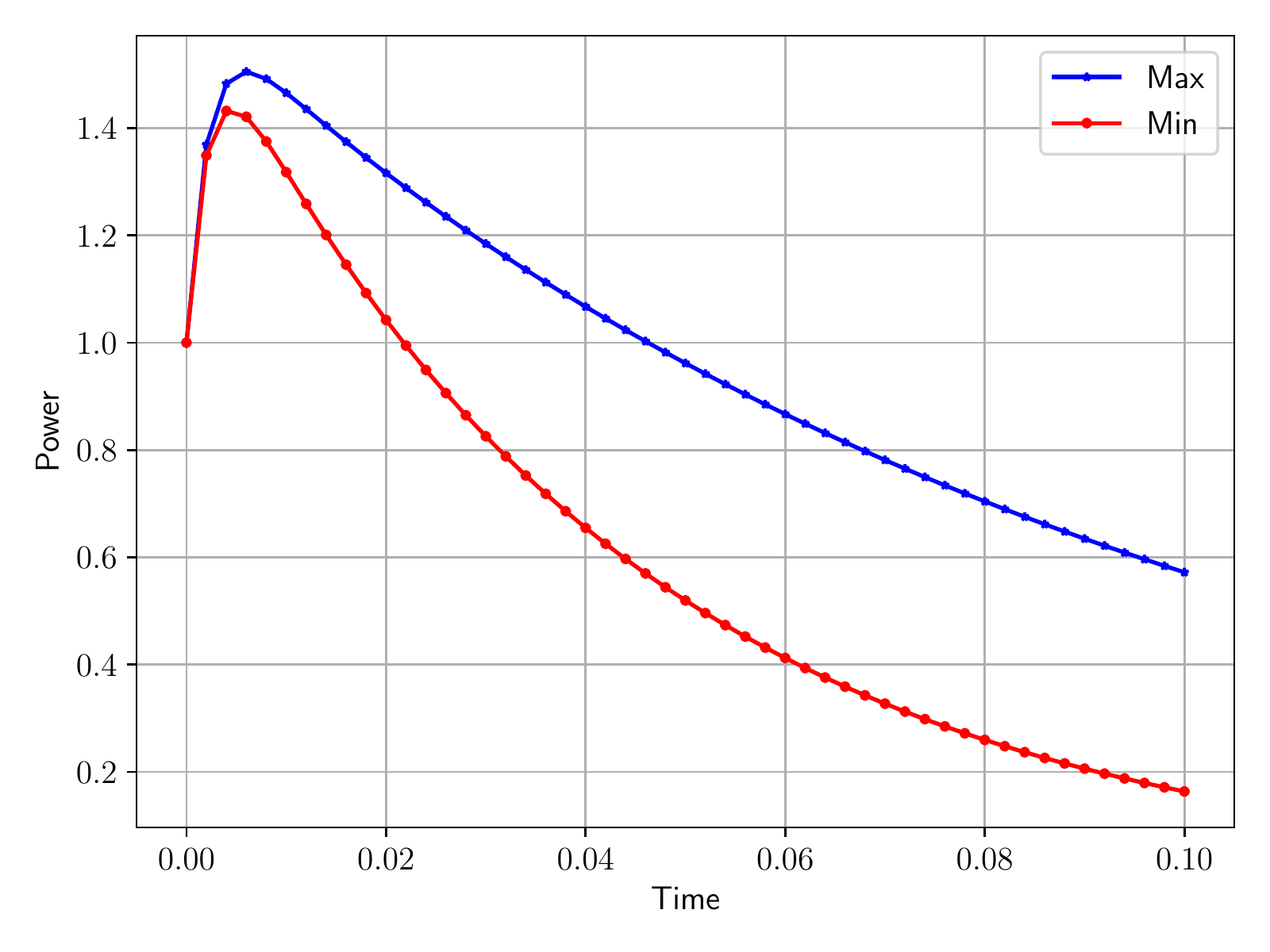}
	\caption{The most (red) and least (blue) subcritical power profiles in the three-dimensional parameter space.}
	\label{fig: Sphere 3D Power Span}
\end{figure}

In this example, the model is parameterized in radius, atom density, and microscopic down-scattering cross-section.
Whereas in the previous section results were presented for simulations which span both subcritical and supercritical regimes, this example is bounded to only subcritical systems.
This is done to be more realistic with respect to the regimes observed in pulsed neutron experiments.
The parameter space is defined by $r_b = [5.94, 6.0]$ cm, $\rho = [0.0495, 0.05]$ atoms/b-cm, and $\sigma_{s, 0 \rightarrow 1} = [1.387, 1.46]$ b and yields the bounding power profiles shown in Figure \ref{fig: Sphere 3D Power Span}.
The relative difference between the minimum and maximum final powers for this parameterization is 111\% and the minimum-to-maximum difference is 250\%.
This space is sampled using a tensor product sampling technique with four equidistant points in each dimension for a total of 64 simulations.

\begin{figure}[!h]
	\centering
	\includegraphics[width=0.7\linewidth]{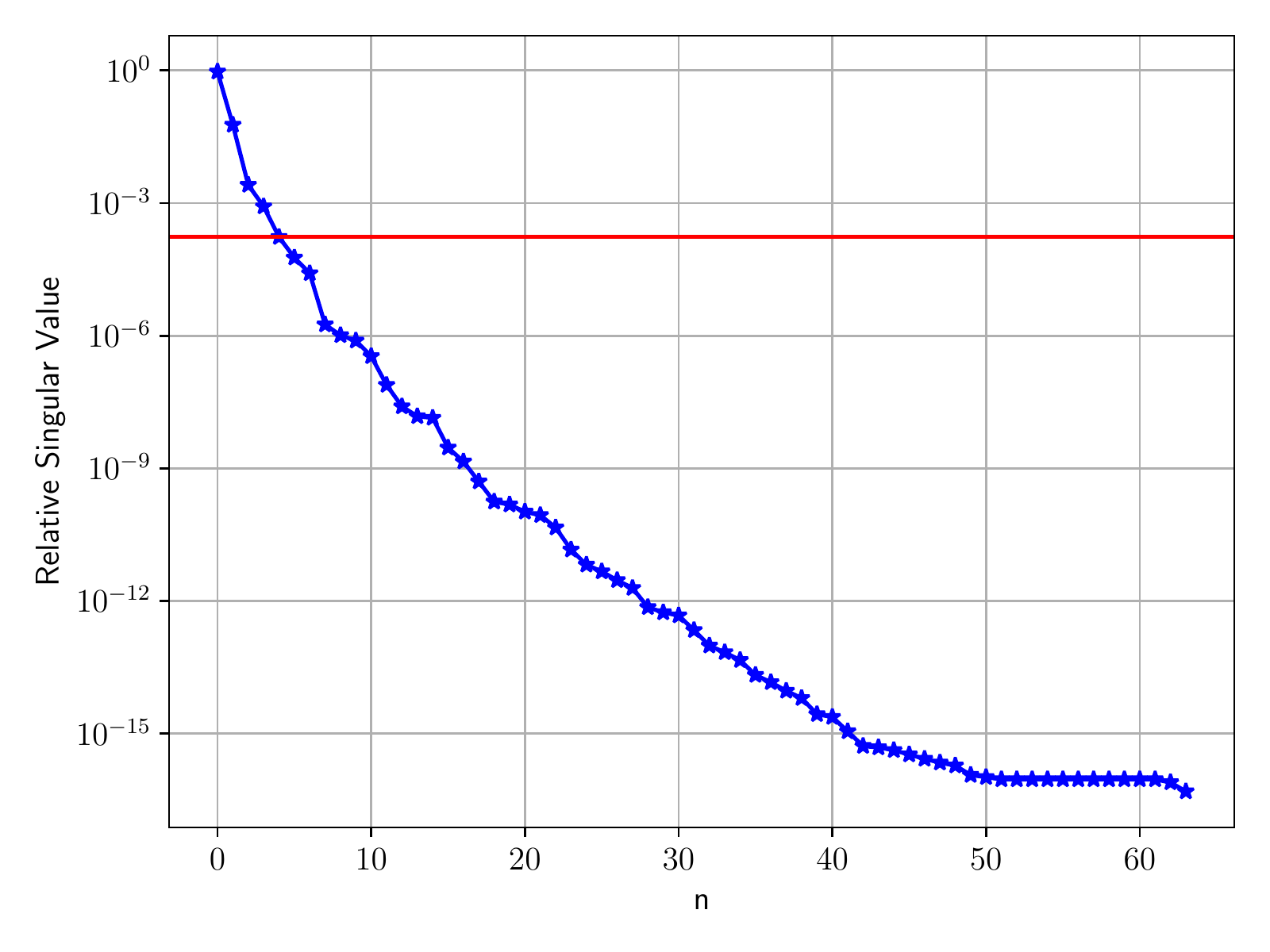}
	\caption{The scree plot for the multi-group scalar flux data.}
	\label{fig: Sphere 3D Multi-Group Scree}
\end{figure}

\begin{figure}[!h]
	\centering
	\includegraphics[width=0.8\linewidth]{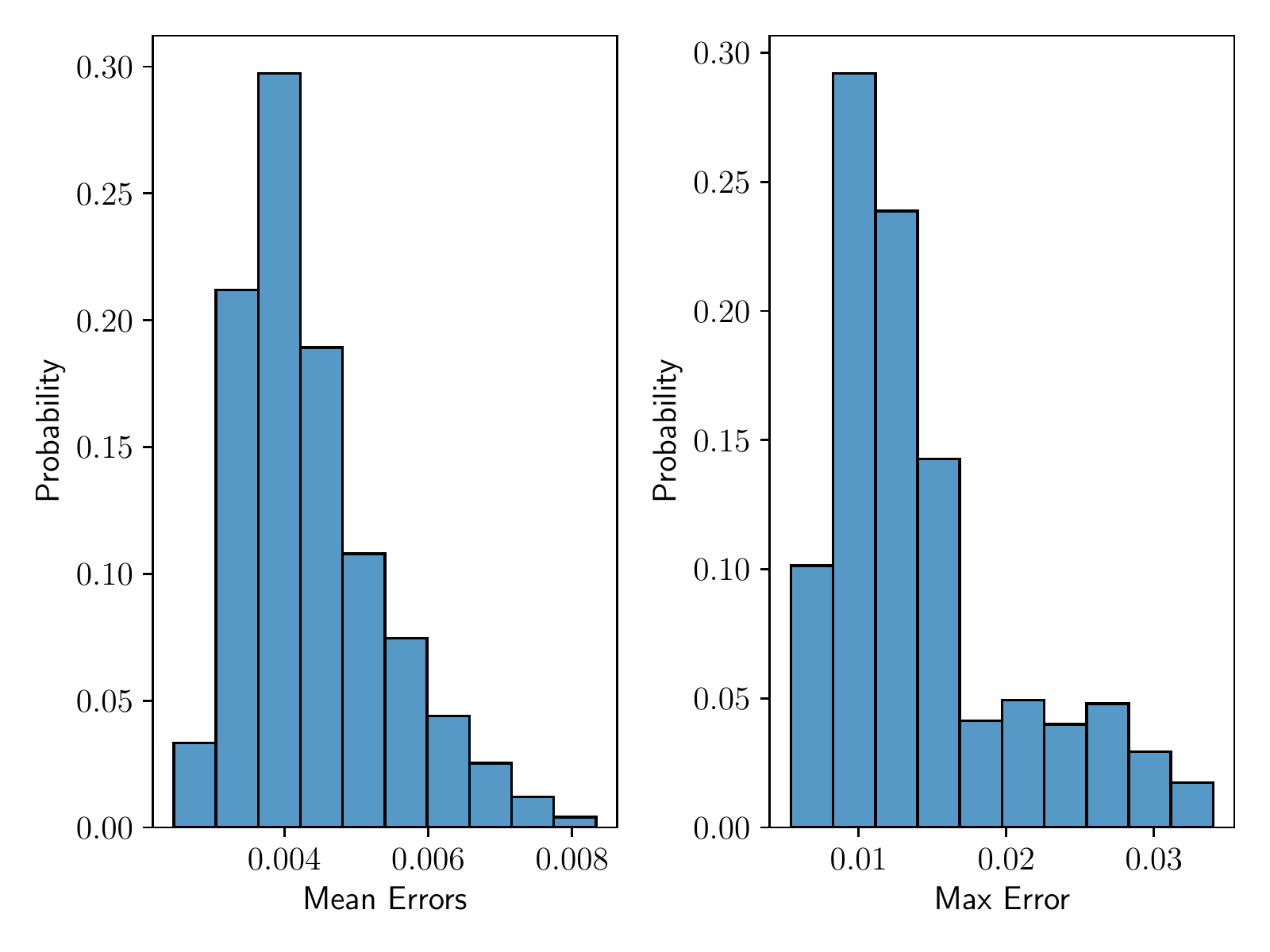}
	\caption{The distribution of mean and maximum errors across 750 cross-validation sets.}
	\label{fig: Sphere 3D Cross Validation}
\end{figure}

The scree-plot for the multi-group scalar flux data is shown in Figure \ref{fig: Sphere 3D Multi-Group Scree}.
In this case, a truncation parameter of $\tau = 10^{-8}$ yields a five mode model.
The reconstruction error over the entire data set is $6.79 \times 10^{-5}$ with mean and maximum simulation reconstruction errors of $5.46 \times 10^{-5}$ and $2.77 \times 10^{-4}$, respectively.
Using leave-one-out cross-validation, the mean prediction error was 0.259\% and the maximum 0.674\%.
This suggests that even with significantly less resolution in the functional behaviors of the mode coefficients, acceptable accuracy can still be attained.
It should be noted that these results are better than those from Section \ref{sec: Sphere 1D} largely due to the restriction to subcritical regime, and therefore, the elimination of the divergent behavior about the subcritical to supercritical transition.
The mean and maximum error distributions across 750 cross-validation sets using $k$-fold cross-validation with $k=3$ are shown in Figure \ref{fig: Sphere 3D Cross Validation}.
Even with a relatively sparsely sampled parameter space and the omission of one-third of those samples from the training set, maximum errors remain under 4\% and average errors remain under 1\%.
As discussed before, by randomly shuffling the available data, one runs the risk of severe under-sampling in some regions of the parameter space. 
In higher dimensions this problem is exacerbated because the space is, in general, sampled much more sparsely due to computational costs.
With this consideration, this result is extremely promising.

\begin{figure}
	\centering
	\includegraphics[width=0.7\linewidth]{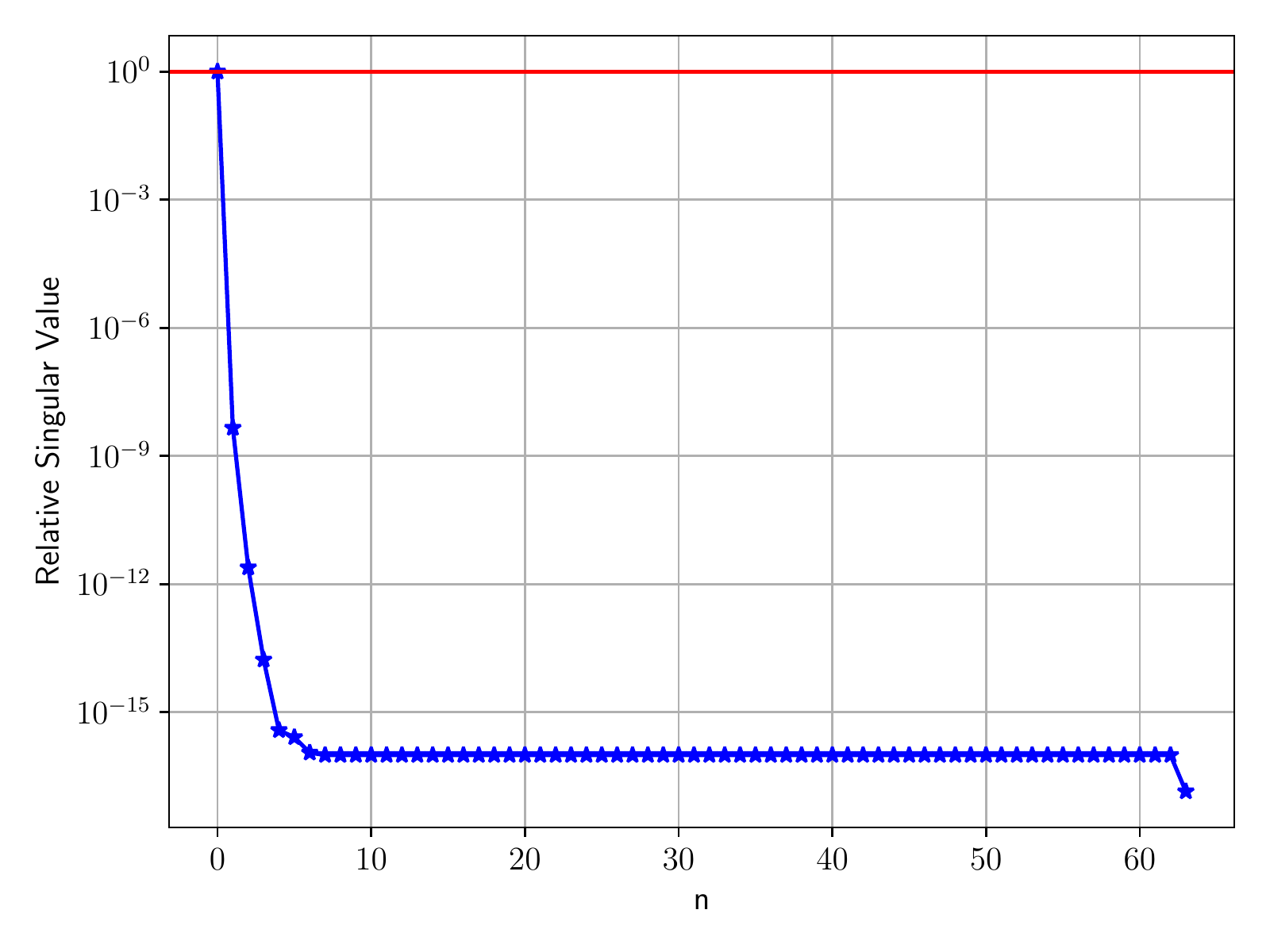}
	\caption{The scree plot for final power density profile data.}
	\label{fig: Sphere 3D Final Power Scree}
\end{figure}

At this stage, the POD-MCI ROM has been shown to accurately predict the parametric behaviors of full-order simulation data.
One of the primary advantages of non-intrusive ROMs, however, is their ability to directly model arbitrary data.
To illustrate this, consider a scenario in which one is interested in only the power density profile at the end of the simulation.
Constructing a ROM from this data reduces the data requirements by a factor equal to the product of the number of groups and time steps.
For this example, this is a factor of 153.
The scree plot for this data is shown in Figure \ref{fig: Sphere 3D Final Power Scree}.
Because at this time in the simulation the solution has settled into its long term behaviors, the power profiles are very similar throughout parameter space and only the magnitude varies.
Because of this, the data can be represented with a truncation limit of $\tau = 10^{-8}$ using only a single POD mode.
Using leave-one-out cross-validation, the mean and maximum errors in predicting the final power density profile are 1.12\% and 3.05\%, respectively, which is still a generally acceptable level of error for a ROM.
The reason for greater error is likely due to the fact that the final power density profile data is ``incomplete'' in that prior dynamic behaviors of system have predictive value for future states.
For this example, the construction time and query time of the ROM both decreased by roughly a factor of four.

As a final example, the power profile as a function of time could be taken as the QoI and directly modeled.
When applying the same procedure, a five mode model is obtained for an energy truncation threshold of $\tau = 10^{-8}$.
Leave-one-out cross-validation yielded results comparable to those when using the full-order solution, but at a four times speedup.
Further, this ROM requires 300 times less data than that using the full-order solution.

\subsection{LRA Benchmark Problem}

\begin{figure}[!h]
	\centering
	\includegraphics[width=0.7\linewidth]{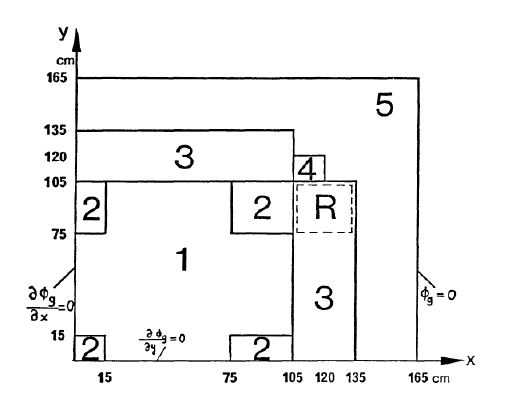}
	\caption{The geometry of the LRA benchmark problem \cite{ANLBenchmarks}.}
	\label{fig: LRA Geometry}
\end{figure}

The LRA benchmark problem is a two-dimension quarter core model of prompt-supercritical excursion with temperature feedback taken from the Argonne National Laboratory Benchmark Problem Book \cite{ANLBenchmarks}.
The geometry for the problem is shown in Figure \ref{fig: LRA Geometry}.
Regions labeled 1-4 are fuel and region 5 is a reflector.
The bottom and left boundaries are reflective and a zero flux condition is imposed on the top and right boundaries. 
This problem uses two energy groups and two delayed neutron precursor species.
The cross-section and delayed neutron precursor data are given in Tables \ref{tab: LRA Cross-Section Data}, \ref{tab: LRA Common Data}, and \ref{tab: LRA Precursor Data}.
\begin{table}[!h]
	\centering
	\caption{Cross-sections for the materials in the LRA benchmark problem.}
	\begin{tabular}{| c | c | c | c | c | c |}
		\hline
		\textbf{Region} & \textbf{Group} & $D$ (cm) & $\sigma_a$ (cm$^{-1}$) & $\nu \sigma_f$ (cm$^{-1}$) & $\sigma_{s, 0 \rightarrow 1}$ (cm$^{-1}$) \\ \hline
		\multirow{2}{*}{1} &  1 & 1.255    & 0.008252 & 0.004602 & \multirow{2}{*}{0.02533} \\		
									   & 2 & 0.211     & 0.1003     & 0.1091       & \\ \hline
		\multirow{2}{*}{2} & 1 & 1.268    & 0.007181  & 0.004609 & \multirow{2}{*}{0.02767} \\		
									   & 2 & 0.1902  & 0.07047   & 0.08675   & \\ \hline		
		\multirow{2}{*}{3} & 1 & 1.259    & 0.008002 & 0.004663 & \multirow{2}{*}{0.2617} \\
									    & 2 & 0.2091 & 0.08344 & 0.1021         & \\ \hline
		\multirow{2}{*}{4} & 1 & 1.259    & 0.008002 & 0.004663 & \multirow{2}{*}{0.2617} \\
										& 2 & 0.2091 & 0.073324 & 0.1021       & \\ \hline		
		\multirow{2}{*}{3} & 1 & 1.257    & 0.006034 & - & \multirow{2}{*}{0.04754} \\
										& 2 & 0.1592 & 0.01911    & -  & \\ \hline			     
	\end{tabular}
	\label{tab: LRA Cross-Section Data}
\end{table}
\begin{table}[!h]
	\centering
	\caption{Nuclear data common to all materials in the LRA benchmark problem.}
	\begin{tabular}{| c | c | c | c | c | c | c | c |}
		\hline
		\textbf{Group} & $v$ (cm/s) & $\chi^p$ & $\chi^d$ & $B_a^2$ (cm$^{-2}$) & $\nu$ & $\nu^p$ & $\nu^d$ \\ \hline
		1 & 3.0e7 & 1 & 1 & \multirow{2}{*}{1.0e-4} & \multirow{2}{*}{2.43} & \multirow{2}{*}{2.41423659} & \multirow{2}{*}{0.01576341} \\
		2& 3.0e5 & 0 & 0 & & & & \\ \hline
	\end{tabular}
	\label{tab: LRA Common Data}
\end{table}
\begin{table}[!h]
	\centering
	\caption{Delayed neutron precursor data for the LRA benchmark problem.}
	\begin{tabular}{| c | c | c |}
		\hline
		\textbf{Species} & $\lambda$ (s$^{-1}$) & Yield Fraction \\ \hline
		1 & 0.0654 & 0.0054 \\ \hline
		2 & 1.35 & 0.001087 \\ \hline
	\end{tabular}
	\label{tab: LRA Precursor Data}
\end{table}
The quantity $B^2_a$ is an axial buckling term used to approximate axial leakage via $B^2_a D_g \phi_g$.

The nominal setup uses 484 7.5 cm $\times$ 7.5 cm cells and is run to 3 s with time steps of 0.01 s. 
The prompt-supercritical excursion is induced via a linear ramp down of the thermal group absorption cross-section in region ``R'' via
\begin{equation}
	\sigma^R_{a,1}(t) = \sigma_{a,1}^3 \begin{cases}
		 1 + \dfrac{t}{t_{\text{ramp}}} \Delta, & 0 < t < t_{\text{ramp}} \\
		1 + \Delta, & t \ge t_{\text{ramp}},
	\end{cases}
\end{equation}
where $\Delta = -0.1212369$ and $t_\text{ramp} = 2.0$ s.
The temperature feedback model is given by
\begin{equation}
	\sigma_{a, 0}(\bm{r}, t) = \left[ 1 + \gamma \left( \sqrt{T(\bm{r}, t)} - \sqrt{T_0(\bm{r}} \right) \right]\sigma^\text{ref}_{a, 0}(\bm{r}),
\end{equation}
were $\gamma = 3.034 \times 10^{-3}$ K$^{1/2}$ and $T_0(\bm{r}) = 300$ K is the initial temperature at position $\bm{r}$.
The nonlinearity is lagged in the simulation.
The temperature is modeled using the adiabatic heat-up model shown in Eq. (\ref{eq: LRA Adiabatic}), where $\alpha = 3.84 \times 10^{-11}$ K-cm$^3$ is a conversion factor which maps the fission rate to a change in temperature.
\begin{equation}
	\ddt{T(\bm{r}, t)} = \alpha \sum_{g=1}^{2} \sigma_{f, g}(\bm{r}) \phi(\bm{r}, t)
	\label{eq: LRA Adiabatic}
\end{equation}
It should be noted that this temperature model implies there is no heat diffusion nor heat removal mechanism.
As a result of this, the fuel monotonically heats up while the reflector remains at the initial temperature throughout the duration of the simulation.
Due to this, no thermal feedback occurs within the reflector.

\begin{figure}[!h]
	\centering
	\begin{subfigure}{0.49\linewidth}
		\centering
		\includegraphics[width=\textwidth]{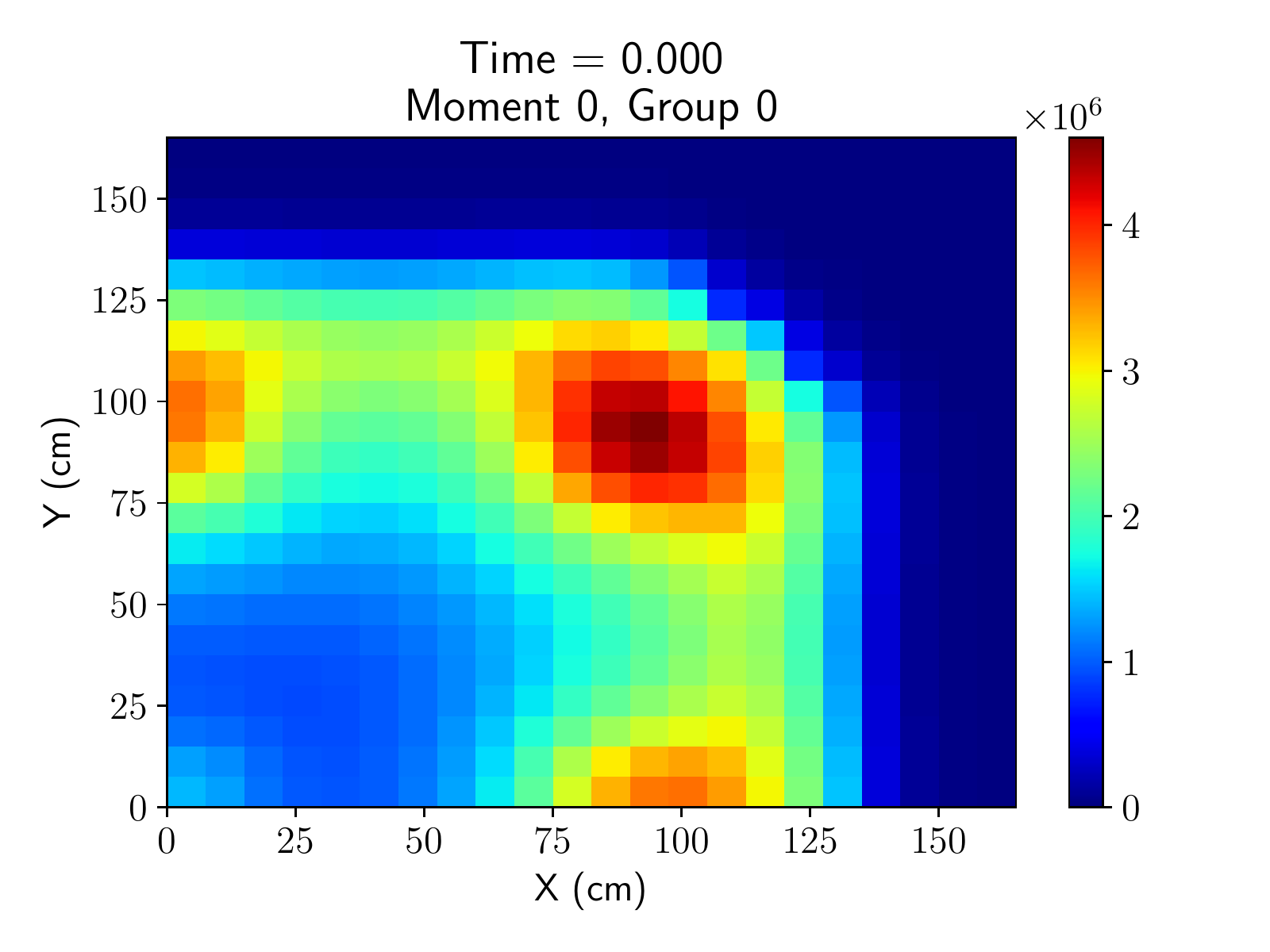}
	\end{subfigure}
	\begin{subfigure}{0.49\linewidth}
		\centering
		\includegraphics[width=\textwidth]{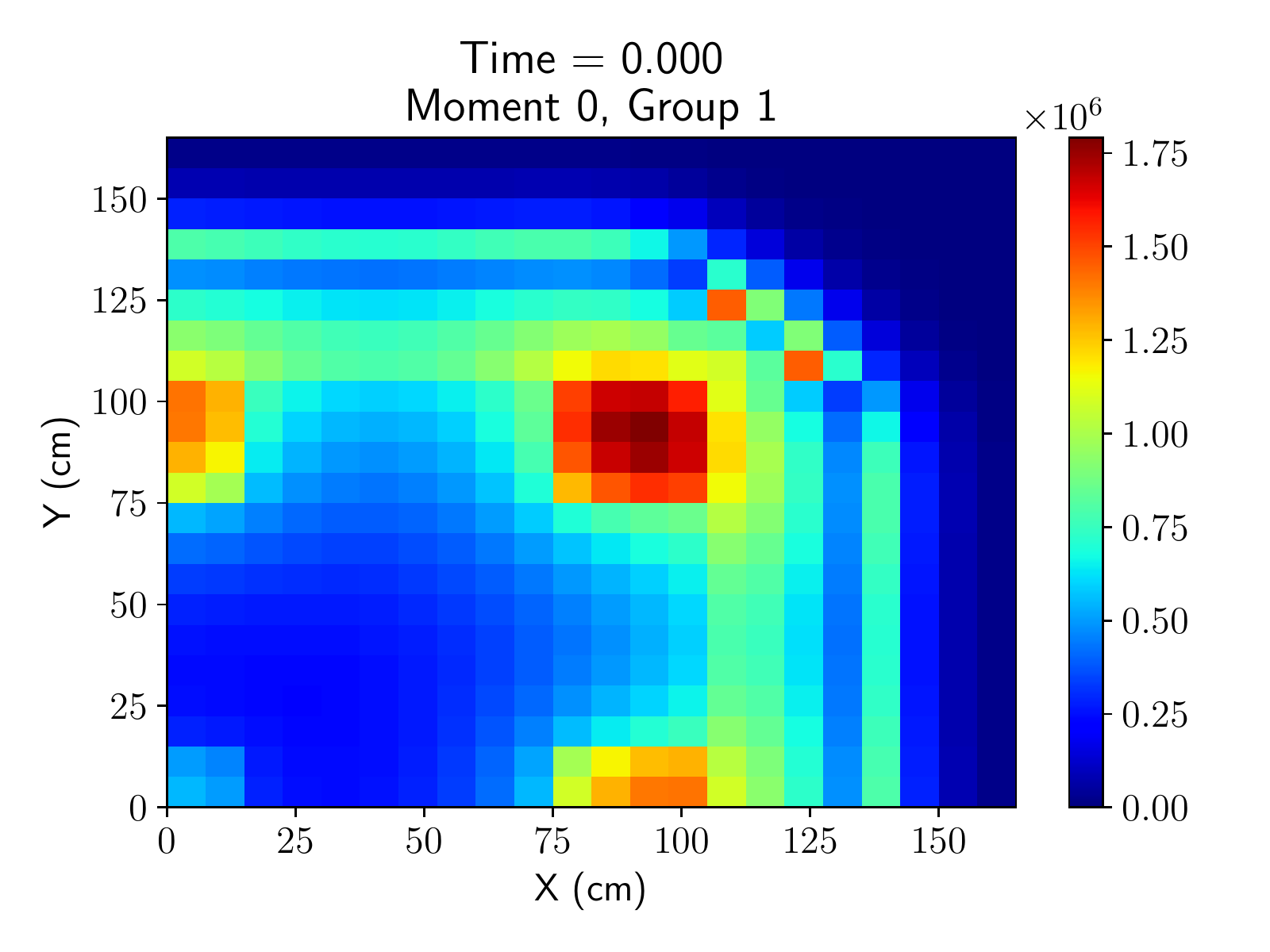}
	\end{subfigure}
	\caption{The initial multi-group scalar flux for the LRA benchmark problem.}
	\label{fig: LRA IC}
\end{figure}

The simulation is initialized by solving the $k$-eigenvalue problem, normalizing the resulting scalar flux distribution such that
\begin{equation}
	\bar{P} = \frac{E_f}{V_\text{core}} \int_{V_\text{core}} \sum_{g=1}^2 \sigma_{f, g}(\bm{r}) \phi_g(\bm{r}, t=0) dV = 10^{-6} \, \text{W cm}^{-3},
\end{equation}
where $\bar{P}$ is the average power density and $E_f = 3.204 \times 10^{-11}$ J is the energy released per fission event.
Further, the fission cross-sections are normalized by the $k$-eigenvalue such that the system is in steady state.
The converged $k$-eigenvalue for this problem was 0.9975.
The initial conditions for each group is shown in Figure \ref{fig: LRA IC}

\begin{figure}[!h]
	\centering
	\begin{subfigure}{0.49\linewidth}
		\centering
		\includegraphics[width=\textwidth]{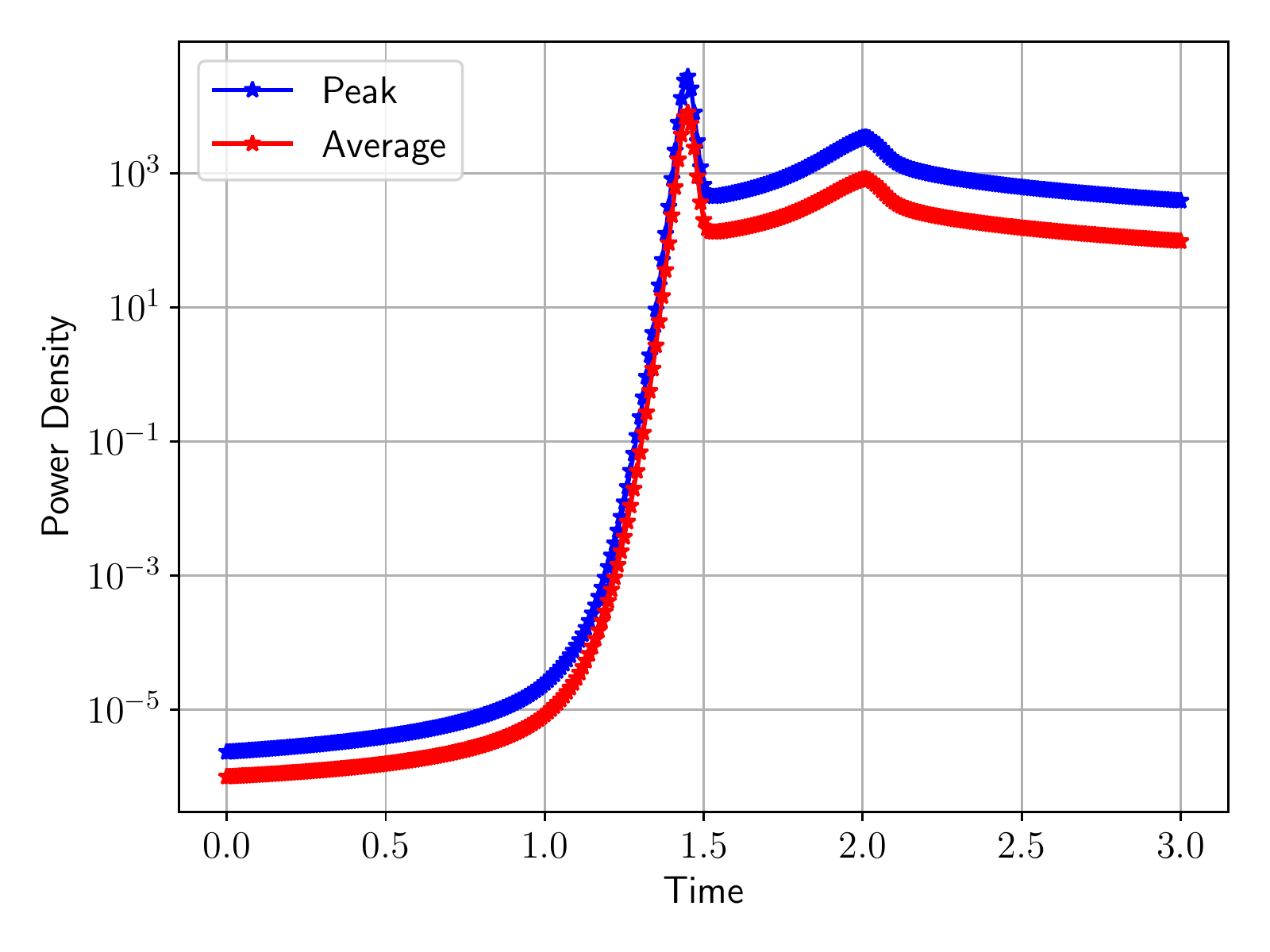}
	\end{subfigure}
	\begin{subfigure}{0.49\linewidth}
		\centering
		\includegraphics[width=\textwidth]{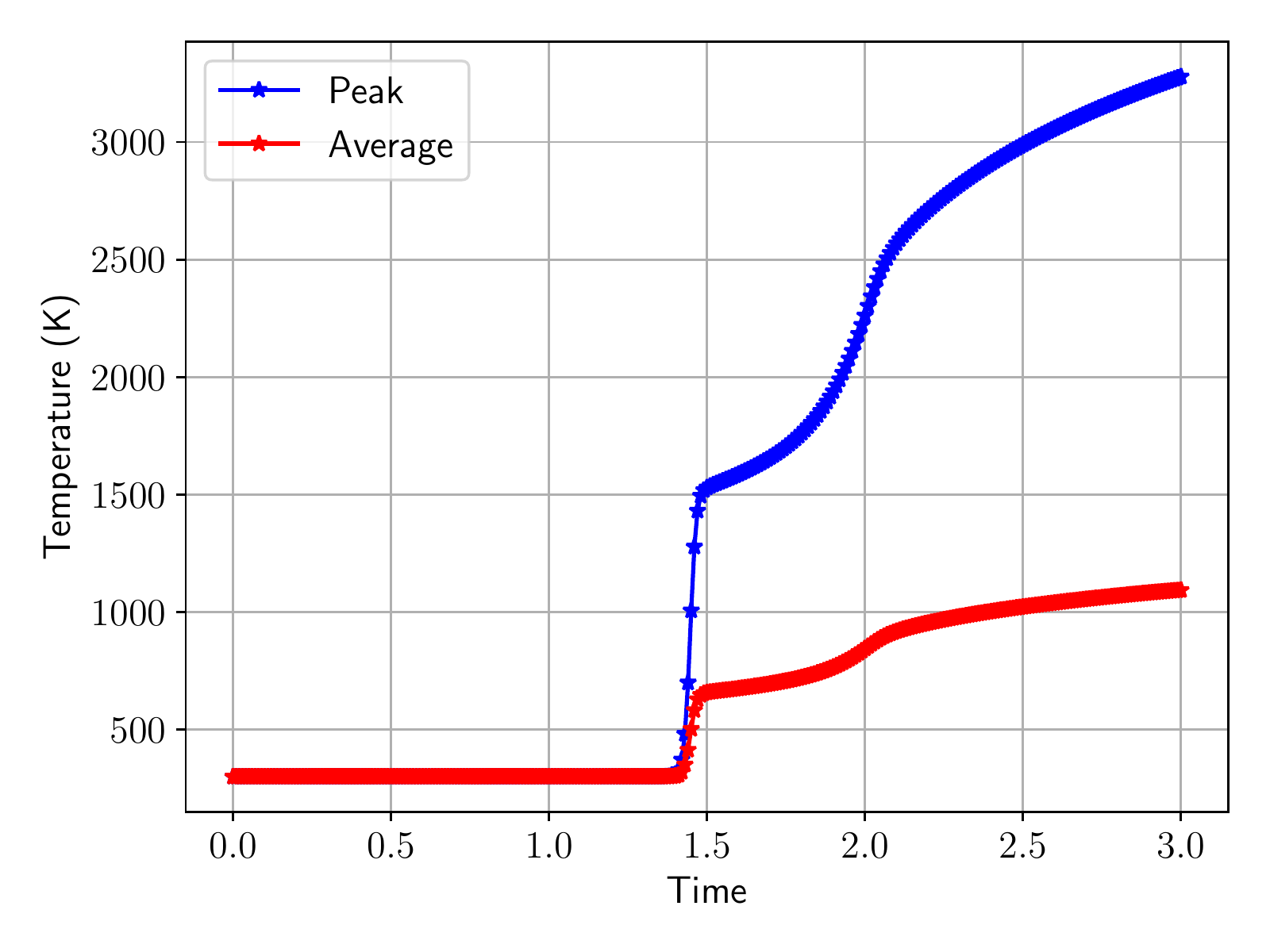}
	\end{subfigure}
	\caption{The average (red) and peak (blue) power density (left) and fuel temperature (right) as a function of time.}
	\label{fig: LRA Dynamics}
\end{figure}

\begin{figure}[!h]
	\centering
	\begin{subfigure}{0.49\linewidth}
		\centering
		\includegraphics[width=\textwidth]{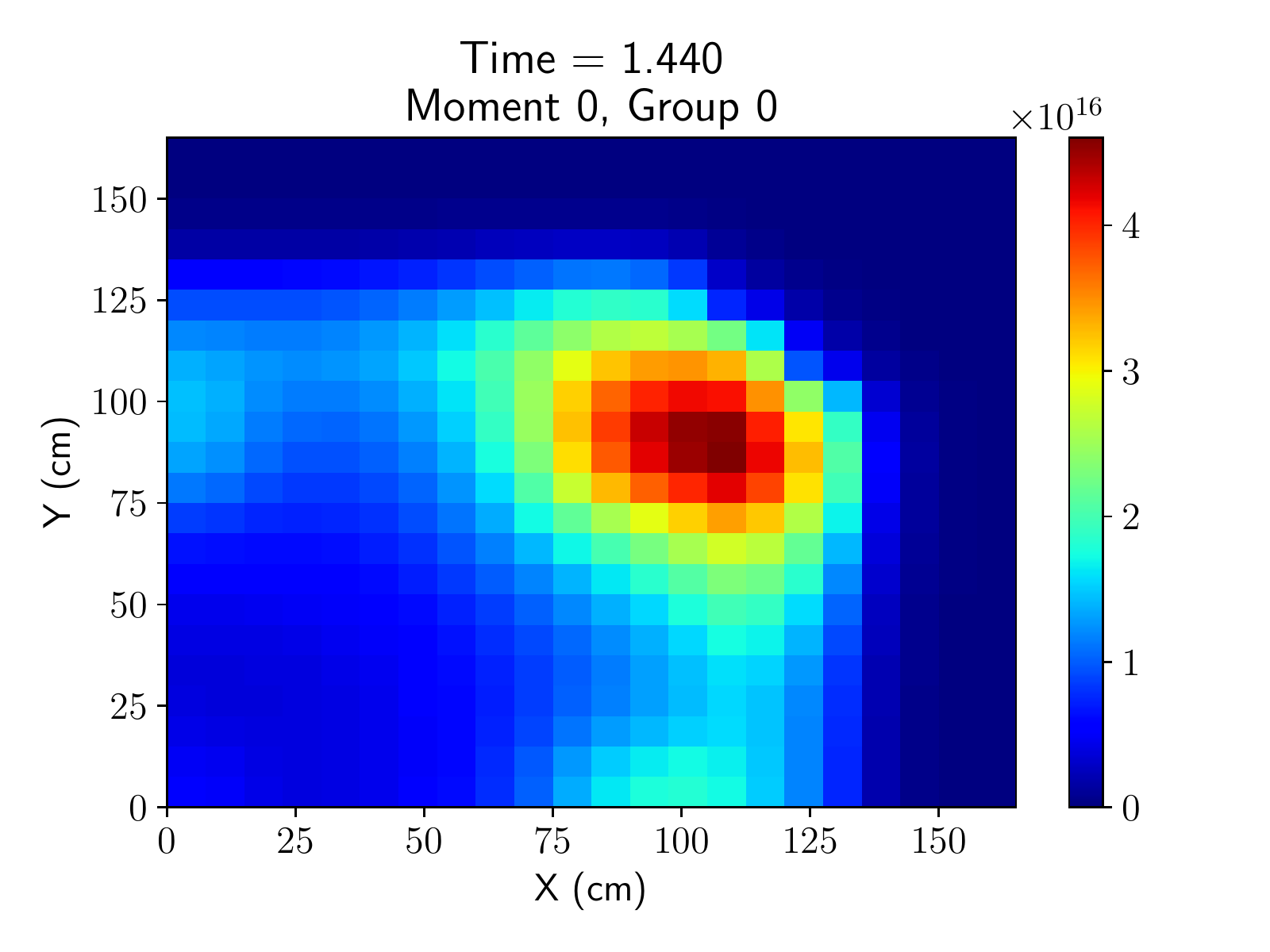}
	\end{subfigure}
	\begin{subfigure}{0.49\linewidth}
		\centering
		\includegraphics[width=\textwidth]{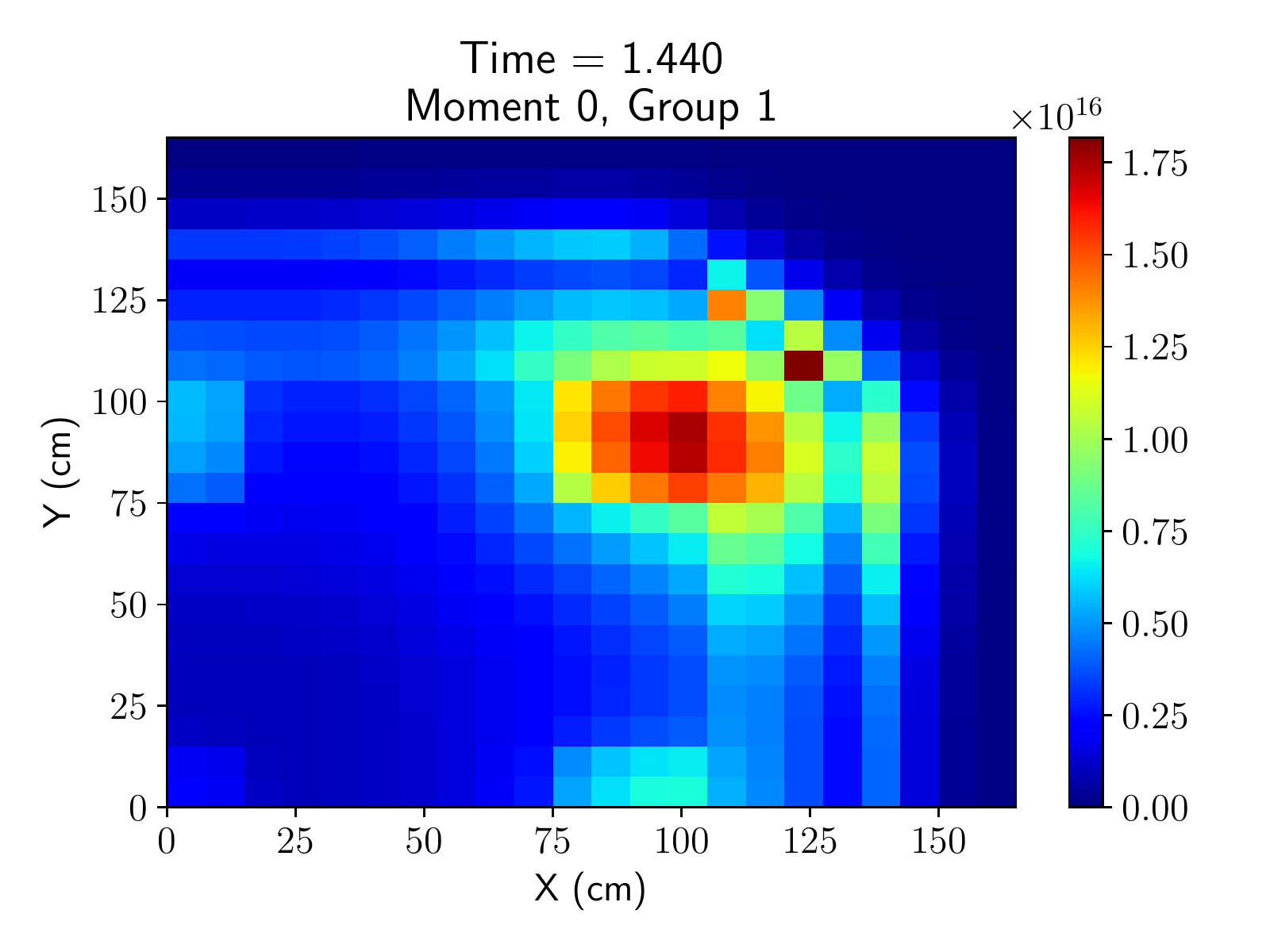}
	\end{subfigure}
	\caption{The multi-group scalar flux at peak power ($t = 1.44$ s) for the LRA benchmark problem.}
	\label{fig: LRA Peak}
\end{figure}

The average and peak power density and fuel temperature as a function of time are shown in Figure \ref{fig: LRA Dynamics}.
For times $t < 1.44$ s, the absorption cross-section in region ``R'' is decreasing, adding reactivity to the system causing the power to increase at increasingly fast rates.
Through this period, the positive reactivity from the decreasing absorption is greater than the negative reactivity contributions from temperature feedback.
At $t = 1.44$ s the power peaks, implying that the temperature feedback has ``caught up'' and is overtaking the positive reactivity from the cross-section ramp.
The multi-group scalar flux profile at peak power time is shown in Figure \ref{fig: LRA Peak}.
For $1.44 < t \le 2$ s, positive reactivity is still being added to the system via the cross-section ramp.
At a time shortly after peak power, the positive reactivity once again overtakes the temperature feedback causing the power trough subsequent rise.
This remains the case until the ramp terminates at $t = 2$ s.
For the remainder of the simulation temperature feedback dominates the problem and the power begins to decay.
The complex dynamics and the ten order of magnitude change in the solution make this benchmark problem an excellent test for the POD-MCI ROM.

\begin{figure}[!h]
	\centering
	\begin{subfigure}{0.49\linewidth}
		\centering
		\includegraphics[width=\textwidth]{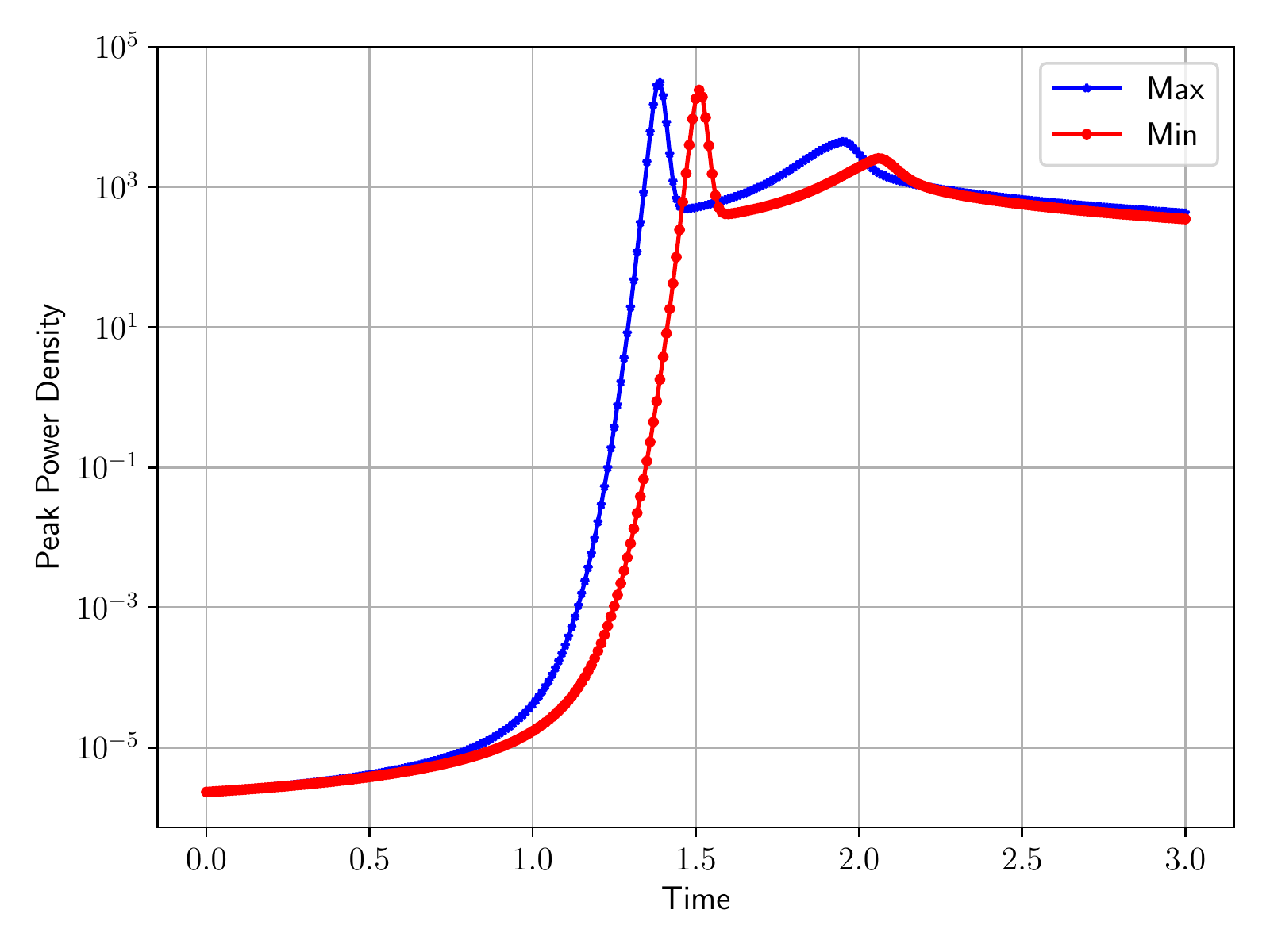}
	\end{subfigure}
	\begin{subfigure}{0.49\linewidth}
		\centering
		\includegraphics[width=\textwidth]{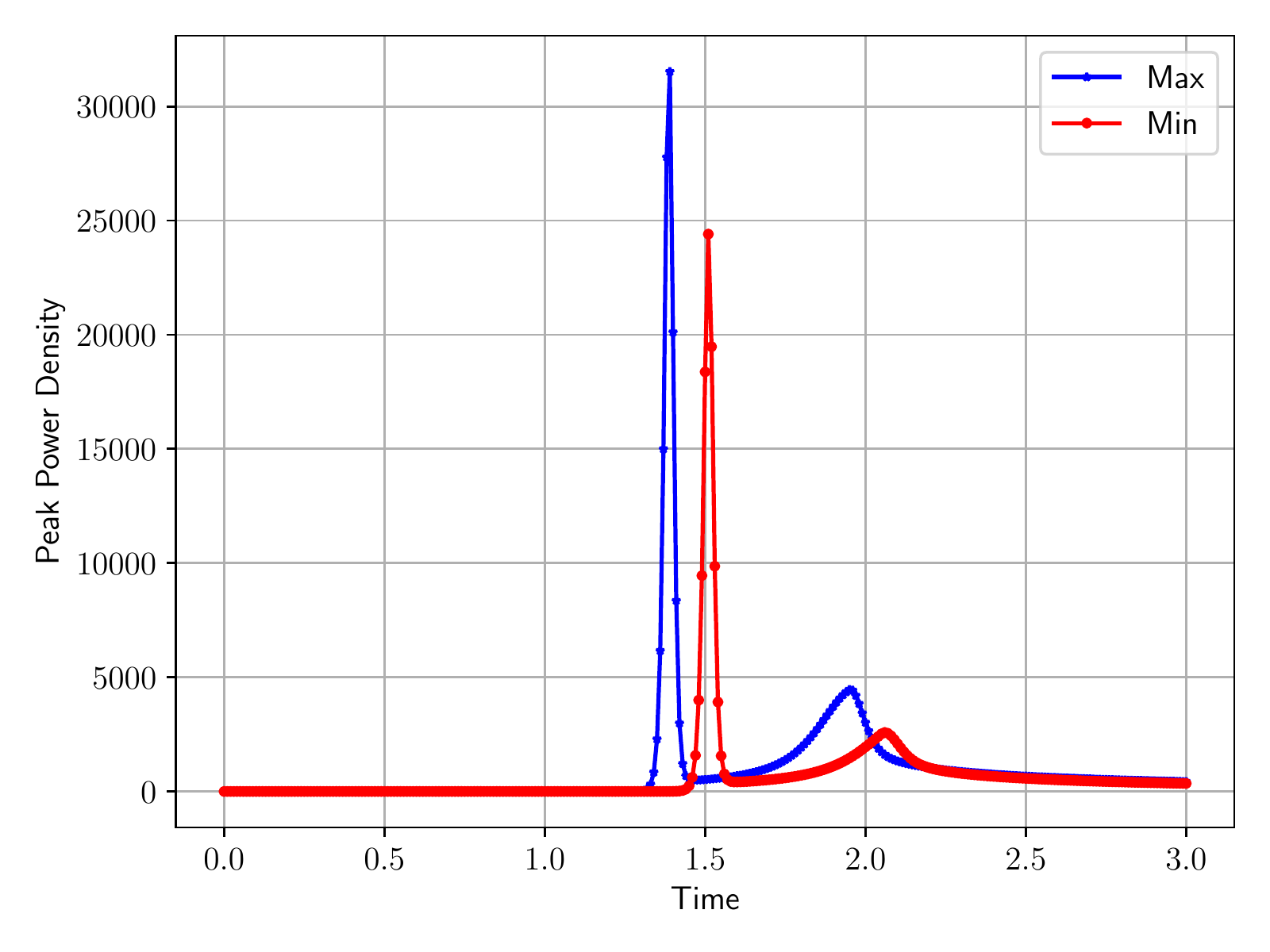}
	\end{subfigure}
	\caption{The power profiles for the simulations with the minimum (red) and maximum (blue) peak power on a log scale (left) and a linear scale (right).}
	\label{fig: LRA Power Span}
\end{figure}

\begin{figure}[!h]
	\centering
	\includegraphics[width=0.7\linewidth]{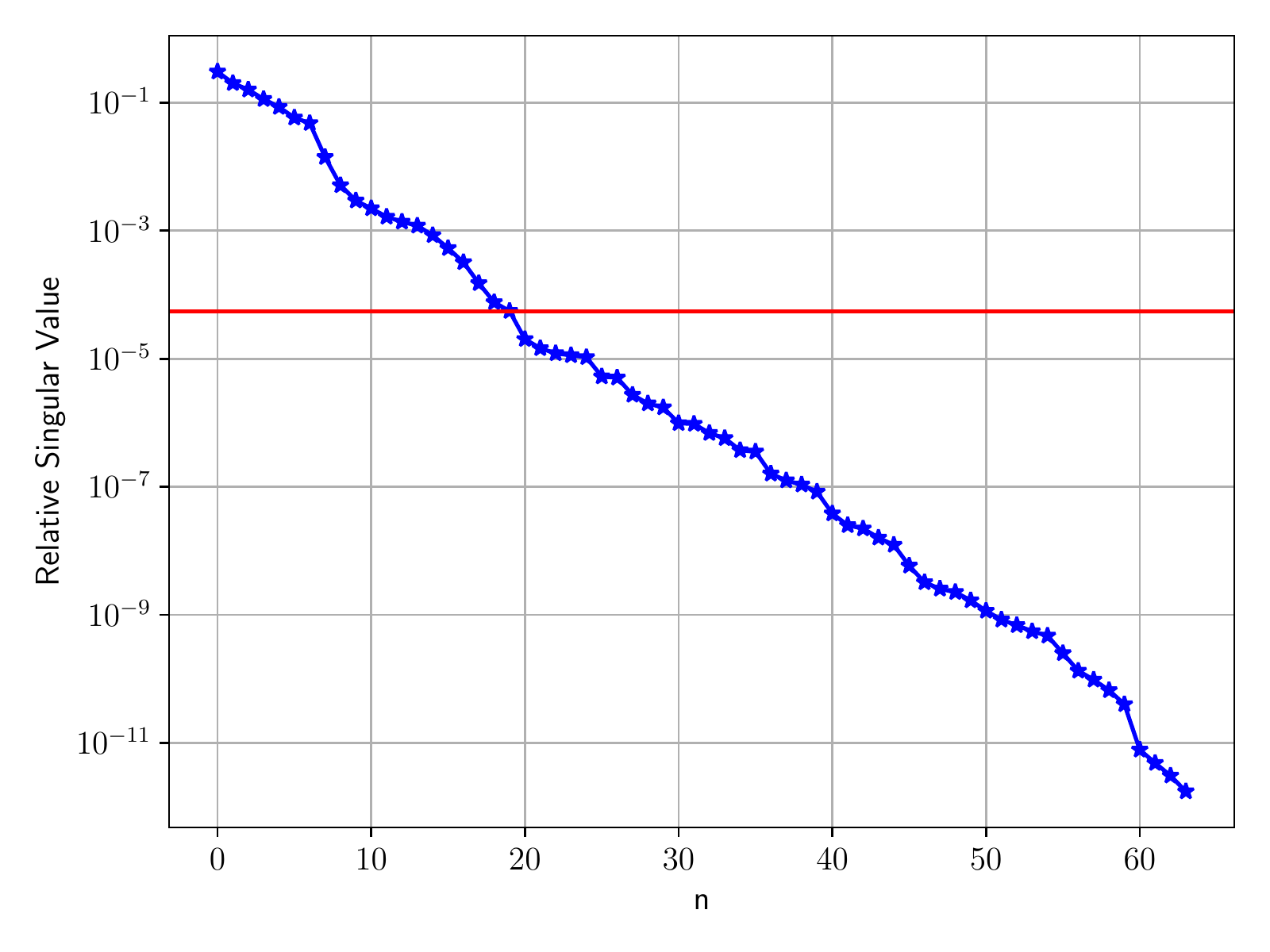}
	\caption{The scree plot for the time-dependent power density profile data for the LRA benchmark problem.}
	\label{fig: LRA Multi-Group SVD}
\end{figure}

In reactor pulse problems, oftentimes one is interested in characterizing the peak power.
To test the POD-MCI ROM, the LRA benchmark problem is parameterized in three variables, the magnitude of the cross-section, $\Delta = -0.1212369 \pm 2.5\%$, the duration of the cross-section ramp $t_\text{ramp} = 2.0 $ s $\pm 2.5\%$, and the feedback coefficient $\gamma = 3.034 \times 10^{-3}$ K$^{1/2}$ $\pm 5\%$.
Each of these parameters directly influence both when the peak power occurs, and its magnitude, in complex ways.
Further, because of the rapid rate of growth leading up to peak power, small perturbations yield orders of magnitude differences in solutions at fixed times prior to the peak power.
This can be seen in Figure \ref{fig: LRA Power Span}.
Because of the wide variation in solution data throughout the parameter space, significantly more POD modes are needed to yield the same amount of energy as in the prior examples.
This is demonstrated for the time-dependent power density profile data in Figure \ref{fig: LRA Multi-Group SVD}.

\begin{figure}[!h]
	\centering
	\includegraphics[width=0.8\linewidth]{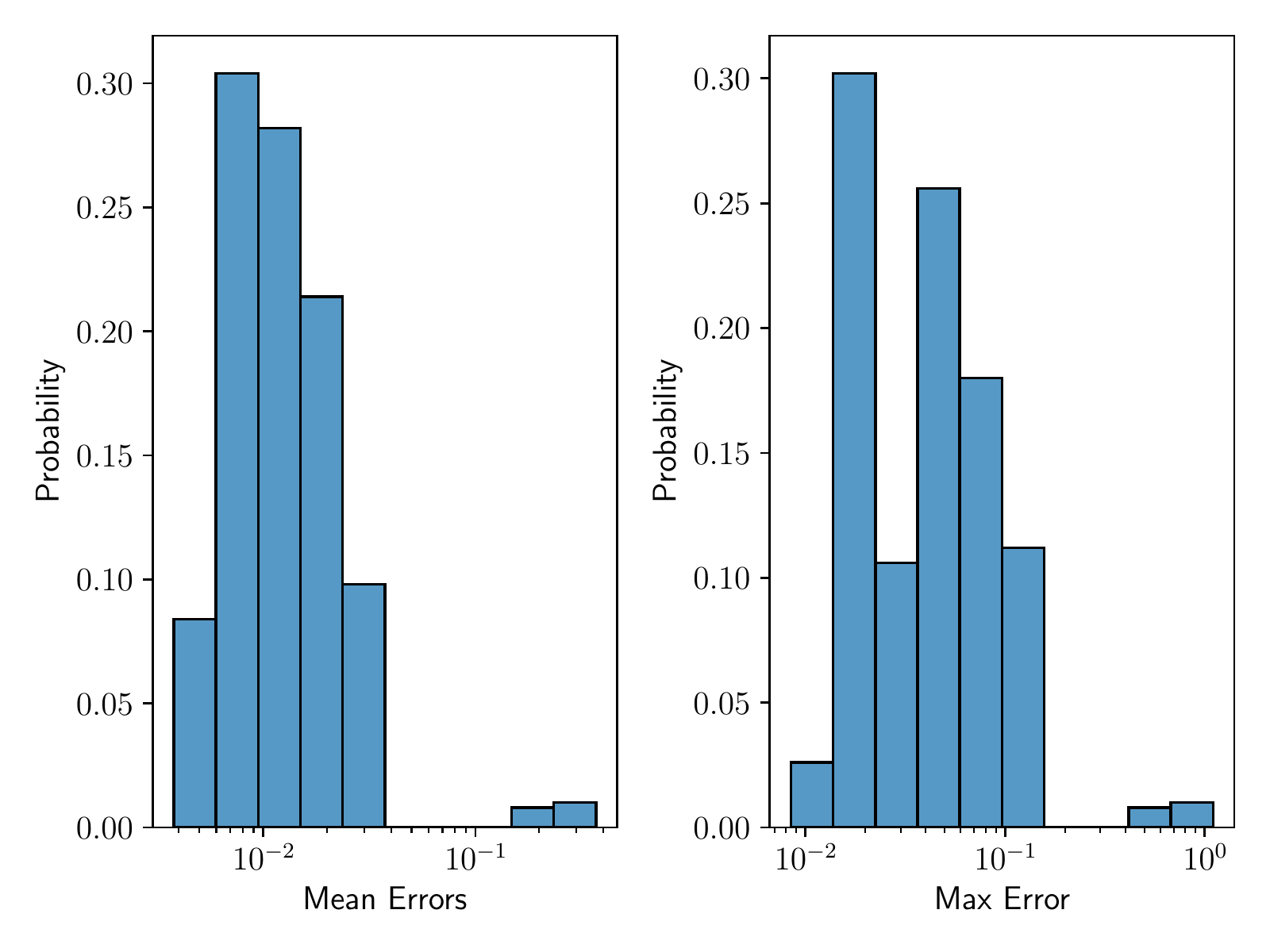}
	\caption{The mean (left) and maximum (right) error distributions across 500 cross-validation sets.}
	\label{fig: LRA Error Distribution}
\end{figure}

While this POD-MCI ROM uses 20 modes for an energy truncation of $\tau = 10^{-8}$, the results are comparable to those from the previous example.
Leave-one-out cross-validation yields mean and maximum prediction errors of 0.657\% and 1.75\%, respectively.
Using three-fold cross-validation as in the previous examples yields poor results due to both the large number of modes and the wide variations in parameter space.
Because the higher-order modes contain more information describing the data set, their contributions to the solution are greater, and therefore, poor interpolations can yield very poor predictions.
With poor sampling in parameter space, as can occur with $k$-fold cross-validation, this problem is exacerbated.
For this reason, five-fold cross-validation is used in lieu of three-fold cross-validation to investigate the effects of the sampling.
This is shown in Figure  \ref{fig: LRA Error Distribution}.
It should be noted that in a very small number of cases, divergent errors are observed.
This is a result of particularly poor sampling of the parameter space within the training set.
Due to the orders of magnitude differences between snapshots at fixed times, poor sampling can yield dramatically larger errors than in the previous example. 
It should be emphasized, however, that maximum observed errors lie between 1.58\% and 12.8\% in 95\% of cases and that large errors can be avoided by performing cross-validation procedures on the available training data before employing the ROM.
In general, these results demonstrate that the performance of the POD-MCI ROM is relatively insensitive to the sampling in parameter space given that an adequate sampling is used.

The FOM for this problem is significantly more expensive than in the previous example.
Without adaptive time stepping or iterating on the nonlinearity, each simulation takes on average 32 s.
With these features, this can more the double.
On average, construction of the POD-MCI ROM took 0.71 s from leave-one-out cross-validation and the average query time was 1.56 ms.
The larger construction time is due to the simulation size.
In this case, there are 301 time steps with 484 grid points, or 145,684 data points per simulation.
It should be noted that this is already an 80\% reduction in the amount of data required for an intrusive ROM where the two-group scalar flux, two precursor species, and temperature would have to be included.
The larger query time is simply a function of the number of modes and the cost of mapping the coefficients back to the full-order result.
With this increased cost, however, the computational speedup factor is approximately 20,500.
A greater speedup factor for a more complex and larger FOM is a favorable property for any ROM.

\begin{figure}[!h]
	\centering
	\includegraphics[width=0.7\linewidth]{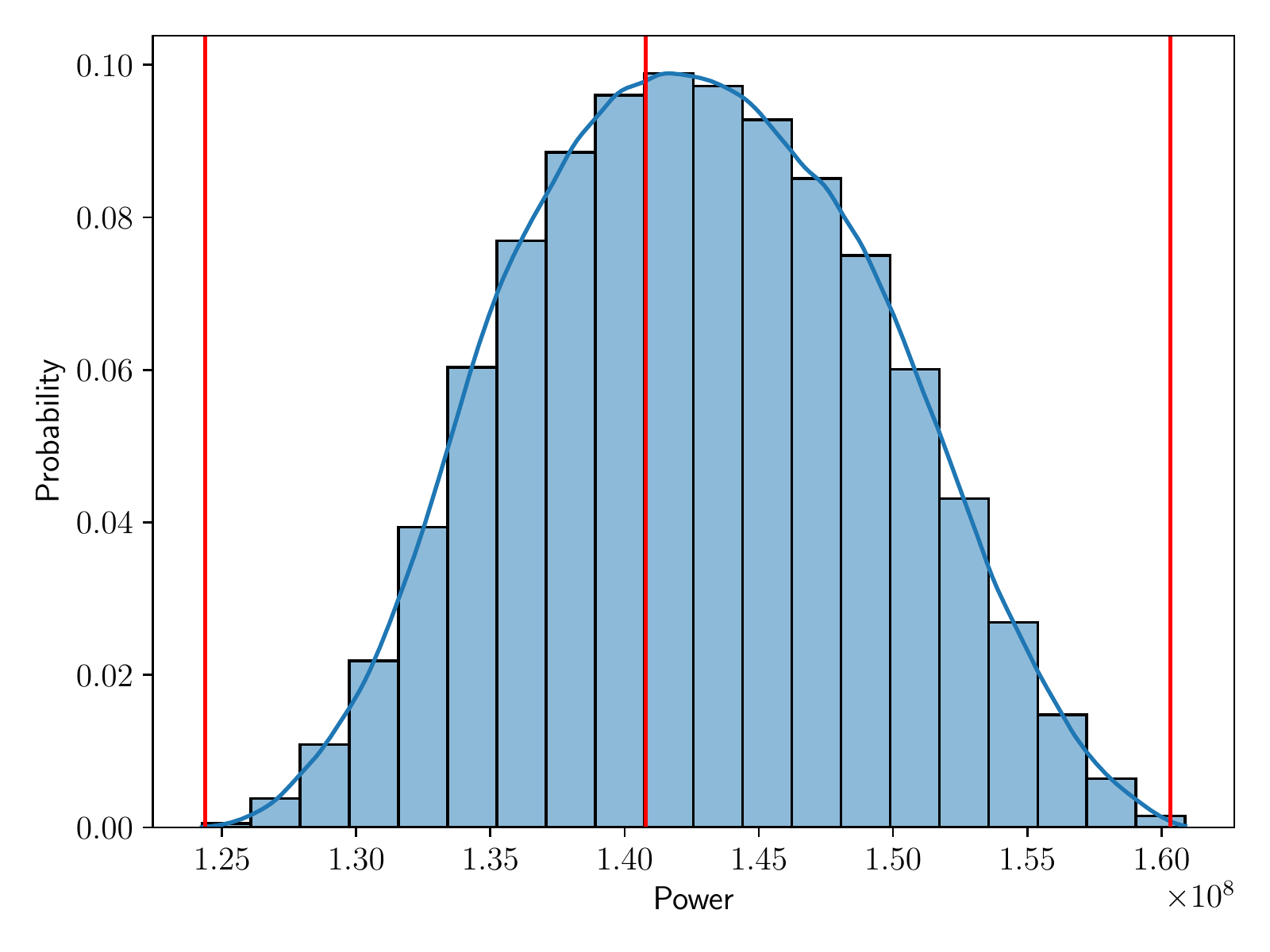}
	\caption{The distribution of peak power across the three-dimensional parameter space plotted with the minimum, maximum, and mean peak power of the available simulation data (red).}
	\label{fig: LRA QoI}
\end{figure}

In a problem such as this, one may want to predict the behaviors of the system at or around peak power.
When using a ROM constructed from the time-dependent power density profile data, the POD-MCI ROM consistently under-predicts the peak power.
The reason for this is relatively simple and is based on the fact that parameter space is treated as a continuum while time is discrete.
Small changes in parameter space likely yields small changes in peak power magnitude and peak power time.
However, peak powers can only be extracted at times corresponding to simulation time steps.
Because of this, if peak power does not lie exactly on a time-step, it will necessarily be under-predicted.
Instead of modeling the time-dependent data, one can simply formulate a ROM based on the peak power density profile for each simulation.
Such a ROM reduces the data requirements by a factor of the number of time steps, or 301, in this case.
As in the previous example, because the underlying behavior of the peak power density profile has lesser variation than the full time-dependent data, the POD model requires significantly fewer modes for the same amount of energy.
For this example, the peak power density profile can be represented with three POD modes and the maximum observed error with leave-one-out cross-validation is 0.424\%.
Taking the space-integrated power density profile as a scalar QoI, the distribution obtained by taking one million random samples in parameter space is shown in Figure \ref{fig: LRA QoI}.
This is clearly a realistic result as that the distribution is approximately normal, lies within the observed range, and is centered near the mean of the observations.
Further, the query time for this ROM was 7.32 $\mu$s across the one million samples, representing a 4.35 million times speedup relative to the FOM.

\section{Conclusions}
\label{sec: Conclusions}

This paper presented the non-intrusive, model agnostic POD-MCI ROM which can be constructed from arbitrary data.
This ROM performs a coordinate transformation from high-dimensional simulation data to a low-rank optimal basis, interpolates in this basis, then maps this result back to the high-dimensional space.
Not only does this method require no modification of FOM source code, but it can greatly alleviate the data requirements of its intrusive counterparts by constructing the ROM only from the simulation data of interest.

A number of different results were presented to characterize the properties and the accuracy of the POD-MCI ROM.
First it was demonstrated the the ROM is stable as a function of mode number.
Because the smoothness of the functional POD mode coefficients degrades with the energy content of the modes, once the energy of the POD modes becomes small relative to the error in interpolating the first few modes, error introduced from poor interpolation becomes negligible.
Second, a bootstrapping methodology was employed to characterize the sensitivity of the method to the parametric sampling.
Outside of cases where regions of parameter space were severely under-sampled, errors remained in an acceptable range.
Lastly, significant improvements in both data requirements and execution time with respect to the FOM were observed when using lower-dimensional data such are power density profiles at particular times or purely time-dependent data.
Each of these results held for both one- and three-dimensional parameter spaces in examples with complex dynamics and significant variations within the parameter space.

This work was largely exploratory and was aimed at characterizing the general performance of the POD-MCI method.
There are many areas of further study that could be explored.
First, optimization on the interpolation technique, such as the RBF kernel function and its hyper-parameters, could yield more accurate results.
Further, techniques such as artificial neural networks or Gaussian processes could be employed to better capture nonlinear behaviors of the mode coefficients.
The latter could also provide useful uncertainty information from the interpolations.
Lastly, exploration into more advanced sampling techniques such as Latin hypercube sampling or sparse grids would provide an indication of the performance of this method on non-regular grids and higher-dimensional spaces.

\pagebreak
\section*{Acknowledgments}

This work was supported by the U.S. Department of Energy  (DOE)  through  the  Los  Alamos  National Laboratory. Los Alamos National Laboratory is operated by Triad National Security, LLC, for the National Nuclear Security  Administration  of  the  DOE  (contract  number 89233218CNA000001). This material is based upon work supported by the National Science Foundation Graduate Research Fellowship under grant number DGE:1746932.Any opinion, findings, and conclusions or recommendations expressed in this material are those of the authors and do not necessarily  reflect  the  views  of  the  National  Science Foundation.

\pagebreak
\bibliographystyle{style/ans_js}                                                                           
\bibliography{bib/background.bib,bib/existing_work.bib}

\begin{thebibliography}{10}
\newcommand{\enquote}[1]{``#1''}
\providecommand{\url}[1]{\texttt{#1}}
\providecommand{\urlprefix}{URL }
\expandafter\ifx\csname urlstyle\endcsname\relax
  \providecommand{\doi}[1]{doi:\discretionary{}{}{}#1}\else
  \providecommand{\doi}{doi:\discretionary{}{}{}\begingroup
  \urlstyle{rm}\Url}\fi

\bibitem{lumley1967structure}
\textsc{J.~L. Lumley}, \enquote{The structure of inhomogeneous turbulent
  flows,} \emph{Atmospheric turbulence and radio wave propagation}, 166--178
  (1967).

\bibitem{sirovich1987turbulenceI}
\textsc{L.~Sirovich}, \enquote{Turbulence and the dynamics of coherent
  structures. I. Coherent structures,} \emph{Quarterly of applied mathematics},
  \textbf{45}, \emph{3}, 561 (1987).

\bibitem{sirovich1987turbulenceII}
\textsc{L.~Sirovich}, \enquote{Turbulence and the dynamics of coherent
  structures. II. Symmetries and transformations,} \emph{Quarterly of Applied
  mathematics}, \textbf{45}, \emph{3}, 573 (1987).

\bibitem{sirovich1987turbulenceIII}
\textsc{L.~Sirovich}, \enquote{Turbulence and the dynamics of coherent
  structures. II. Symmetries and transformations,} \emph{Quarterly of Applied
  mathematics}, \textbf{45}, \emph{3}, 573 (1987).

\bibitem{buchan2013pod}
\textsc{A.~Buchan}, \textsc{C.~Pain}, \textsc{F.~Fang}, and \textsc{I.~Navon},
  \enquote{A POD reduced-order model for eigenvalue problems with application
  to reactor physics,} \emph{International Journal for Numerical Methods in
  Engineering}, \textbf{95}, \emph{12}, 1011 (2013).

\bibitem{sartori2016multi}
\textsc{A.~Sartori}, \textsc{A.~Cammi}, \textsc{L.~Luzzi}, and
  \textsc{G.~Rozza}, \enquote{A multi-physics reduced order model for the
  analysis of Lead Fast Reactor single channel,} \emph{Annals of Nuclear
  Energy}, \textbf{87}, 198 (2016).

\bibitem{german2019reduced}
\textsc{P.~German} and \textsc{J.~C. Ragusa}, \enquote{Reduced-order modeling
  of parameterized multi-group diffusion k-eigenvalue problems,} \emph{Annals
  of Nuclear Energy}, \textbf{134}, 144 (2019).

\bibitem{german2019application}
\textsc{P.~German}, \textsc{J.~C. Ragusa}, and \textsc{C.~Fiorina},
  \enquote{Application of multiphysics model order reduction to
  doppler/neutronic feedback,} \emph{EPJ Nuclear Sciences \& Technologies},
  \textbf{5}, \emph{ARTICLE}, 17 (2019).

\bibitem{behne2022minimally}
\textsc{P.~Behne}, \textsc{J.~Vermaak}, and \textsc{J.~C. Ragusa},
  \enquote{Minimally-Invasive Parametric Model-Order Reduction for Sweep-Based
  Radiation Transport,} \emph{Journal of Computational Physics}, 111525 (2022).

\bibitem{prill2014semi}
\textsc{D.~Prill} and \textsc{A.~Class}, \enquote{Semi-automated proper
  orthogonal decomposition reduced order model non-linear analysis for future
  BWR stability,} \emph{Annals of Nuclear Energy}, \textbf{67}, 70 (2014).

\bibitem{sartori2014comparison}
\textsc{A.~Sartori}, \textsc{D.~Baroli}, \textsc{A.~Cammi}, \textsc{D.~Chiesa},
  \textsc{L.~Luzzi}, \textsc{R.~Ponciroli}, \textsc{E.~Previtali},
  \textsc{M.~E. Ricotti}, \textsc{G.~Rozza}, and \textsc{M.~Sisti},
  \enquote{Comparison of a Modal Method and a Proper Orthogonal Decomposition
  approach for multi-group time-dependent reactor spatial kinetics,}
  \emph{Annals of Nuclear Energy}, \textbf{71}, 217 (2014).

\bibitem{soucasse2019angular}
\textsc{L.~Soucasse}, \textsc{A.~G. Buchan}, \textsc{S.~Dargaville}, and
  \textsc{C.~C. Pain}, \enquote{An angular reduced order model for radiative
  transfer in non grey media,} \emph{Journal of Quantitative Spectroscopy and
  Radiative Transfer}, \textbf{229}, 23 (2019).

\bibitem{sun2020pod}
\textsc{Y.~Sun}, \textsc{J.~Yang}, \textsc{Y.~Wang}, \textsc{Z.~Li}, and
  \textsc{Y.~Ma}, \enquote{A POD reduced-order model for resolving the neutron
  transport problems of nuclear reactor,} \emph{Annals of Nuclear Energy},
  \textbf{149}, 107799 (2020).

\bibitem{hughes2020discontinuous}
\textsc{A.~C. Hughes} and \textsc{A.~G. Buchan}, \enquote{A discontinuous and
  adaptive reduced order model for the angular discretization of the Boltzmann
  transport equation,} \emph{International Journal for Numerical Methods in
  Engineering}, \textbf{121}, \emph{24}, 5647 (2020).

\bibitem{tano2021affine}
\textsc{M.~Tano}, \textsc{J.~Ragusa}, \textsc{D.~Caron}, and \textsc{P.~Behne},
  \enquote{Affine reduced-order model for radiation transport problems in
  cylindrical coordinates,} \emph{Annals of Nuclear Energy}, \textbf{158},
  108214 (2021).

\bibitem{permann2020moose}
\textsc{C.~J. Permann}, \textsc{D.~R. Gaston}, \textsc{D.~Andr{\v{s}}},
  \textsc{R.~W. Carlsen}, \textsc{F.~Kong}, \textsc{A.~D. Lindsay},
  \textsc{J.~M. Miller}, \textsc{J.~W. Peterson}, \textsc{A.~E. Slaughter},
  \textsc{R.~H. Stogner} \textsc{et~al.}, \enquote{MOOSE: Enabling massively
  parallel multiphysics simulation,} \emph{SoftwareX}, \textbf{11}, 100430
  (2020).

\bibitem{fiorina2015gen}
\textsc{C.~Fiorina}, \textsc{I.~Clifford}, \textsc{M.~Aufiero}, and
  \textsc{K.~Mikityuk}, \enquote{GeN-Foam: a novel OpenFOAM{\textregistered}
  based multi-physics solver for 2D/3D transient analysis of nuclear reactors,}
  \emph{Nuclear Engineering and Design}, \textbf{294}, 24 (2015).

\bibitem{chaturantabut2010nonlinear}
\textsc{S.~Chaturantabut} and \textsc{D.~C. Sorensen}, \enquote{Nonlinear model
  reduction via discrete empirical interpolation,} \emph{SIAM Journal on
  Scientific Computing}, \textbf{32}, \emph{5}, 2737 (2010).

\bibitem{khuwaileh2020gaussian}
\textsc{B.~A. Khuwaileh} and \textsc{W.~A. Metwally}, \enquote{Gaussian process
  approach for dose mapping in radiation fields,} \emph{Nuclear Engineering and
  Technology}, \textbf{52}, \emph{8}, 1807 (2020).

\bibitem{baraldi2015prognostics}
\textsc{P.~Baraldi}, \textsc{F.~Mangili}, and \textsc{E.~Zio}, \enquote{A
  prognostics approach to nuclear component degradation modeling based on
  Gaussian Process Regression,} \emph{Progress in Nuclear Energy}, \textbf{78},
  141 (2015).

\bibitem{wu2016inverse}
\textsc{X.~Wu} and \textsc{T.~Kozlowski}, \enquote{Inverse uncertainty
  quantification of reactor simulation with polynomial chaos surrogate model,}
  \emph{Transactions of American Nuclear Society. American Nuclear Society, New
  Orleans, LA, USA} (2016).

\bibitem{hardy2019dynamic}
\textsc{Z.~K. Hardy}, \textsc{J.~E. Morel}, and \textsc{C.~Ahrens},
  \enquote{Dynamic mode decomposition for subcritical metal systems,}
  \emph{Nuclear Science and Engineering}, \textbf{193}, \emph{11}, 1173 (2019).

\bibitem{mcclarren2019calculating}
\textsc{R.~G. McClarren}, \enquote{Calculating time eigenvalues of the neutron
  transport equation with dynamic mode decomposition,} \emph{Nuclear Science
  and Engineering}, \textbf{193}, \emph{8}, 854 (2019).

\bibitem{di2020dynamic}
\textsc{A.~Di~Ronco}, \textsc{C.~Introini}, \textsc{E.~Cervi},
  \textsc{S.~Lorenzi}, \textsc{Y.~S. Jeong}, \textsc{S.~B. Seo}, \textsc{I.~C.
  Bang}, \textsc{F.~Giacobbo}, and \textsc{A.~Cammi}, \enquote{Dynamic mode
  decomposition for the stability analysis of the Molten Salt Fast Reactor
  core,} \emph{Nuclear Engineering and Design}, \textbf{362}, 110529 (2020).

\bibitem{abdo2019modeling}
\textsc{M.~Abdo}, \textsc{R.~Elzohery}, and \textsc{J.~A. Roberts},
  \enquote{Modeling isotopic evolution with surrogates based on dynamic mode
  decomposition,} \emph{Annals of Nuclear Energy}, \textbf{129}, 280 (2019).

\bibitem{abdo2018analysis}
\textsc{M.~Abdo}, \textsc{R.~Elzohery}, and \textsc{J.~A. Roberts},
  \enquote{Analysis of the LRA Reactor Benchmark Using Dynamic Mode
  Decomposition,} \emph{Trans. Am. Nucl. Soc}, \textbf{119}, 683 (2018).

\bibitem{mckay2000comparison}
\textsc{M.~D. McKay}, \textsc{R.~J. Beckman}, and \textsc{W.~J. Conover},
  \enquote{A comparison of three methods for selecting values of input
  variables in the analysis of output from a computer code,}
  \emph{Technometrics}, \textbf{42}, \emph{1}, 55 (2000).

\bibitem{smolyak1963quadrature}
\textsc{S.~A. Smolyak}, \enquote{Quadrature and interpolation formulas for
  tensor products of certain classes of functions,} \emph{Doklady Akademii
  Nauk}, vol. 148, 1042--1045, Russian Academy of Sciences (1963).

\bibitem{rathinam2003new}
\textsc{M.~Rathinam} and \textsc{L.~R. Petzold}, \enquote{A new look at proper
  orthogonal decomposition,} \emph{SIAM Journal on Numerical Analysis},
  \textbf{41}, \emph{5}, 1893 (2003).

\bibitem{edwards2011nonlinear}
\textsc{J.~D. Edwards}, \textsc{J.~E. Morel}, and \textsc{D.~A. Knoll},
  \enquote{Nonlinear variants of the TR/BDF2 method for thermal radiative
  diffusion,} \emph{Journal of Computational Physics}, \textbf{230}, \emph{4},
  1198 (2011).

\bibitem{SciPy}
\textsc{P.~Virtanen}, \textsc{R.~Gommers}, \textsc{T.~E. Oliphant},
  \textsc{M.~Haberland}, \textsc{T.~Reddy}, \textsc{D.~Cournapeau},
  \textsc{E.~Burovski}, \textsc{P.~Peterson}, \textsc{W.~Weckesser},
  \textsc{J.~Bright}, \textsc{S.~J. {van der Walt}}, \textsc{M.~Brett},
  \textsc{J.~Wilson}, \textsc{K.~J. Millman}, \textsc{N.~Mayorov},
  \textsc{A.~R.~J. Nelson}, \textsc{E.~Jones}, \textsc{R.~Kern},
  \textsc{E.~Larson}, \textsc{C.~J. Carey}, \textsc{{\.I}.~Polat},
  \textsc{Y.~Feng}, \textsc{E.~W. Moore}, \textsc{J.~{VanderPlas}},
  \textsc{D.~Laxalde}, \textsc{J.~Perktold}, \textsc{R.~Cimrman},
  \textsc{I.~Henriksen}, \textsc{E.~A. Quintero}, \textsc{C.~R. Harris},
  \textsc{A.~M. Archibald}, \textsc{A.~H. Ribeiro}, \textsc{F.~Pedregosa},
  \textsc{P.~{van Mulbregt}}, and \textsc{{SciPy 1.0 Contributors}},
  \enquote{{{SciPy} 1.0: Fundamental Algorithms for Scientific Computing in
  Python},} \emph{Nature Methods}, \textbf{17}, 261 (2020);
  {10.1038/s41592-019-0686-2}.

\bibitem{Scikit-Learn}
\textsc{F.~Pedregosa}, \textsc{G.~Varoquaux}, \textsc{A.~Gramfort},
  \textsc{V.~Michel}, \textsc{B.~Thirion}, \textsc{O.~Grisel},
  \textsc{M.~Blondel}, \textsc{P.~Prettenhofer}, \textsc{R.~Weiss},
  \textsc{V.~Dubourg}, \textsc{J.~Vanderplas}, \textsc{A.~Passos},
  \textsc{D.~Cournapeau}, \textsc{M.~Brucher}, \textsc{M.~Perrot}, and
  \textsc{E.~Duchesnay}, \enquote{Scikit-learn: Machine Learning in {P}ython,}
  \emph{Journal of Machine Learning Research}, \textbf{12}, 2825 (2011).

\bibitem{ANLBenchmarks}
\textsc{A.~C. Center}, \enquote{Benchmark Problem Book, ANL7416-Suppl. 1, 2,}
  (1977).

\end{thebibliography}

\end{document}